\def\yes{\if00}
\def\no{\if01}
\def\iftwelvept{\yes}
\def\ifusepdf{\no}
\def\ifusepsfont{\no}

\iftwelvept
\documentclass[leqno,12pt]{amsart}
\else
\documentclass[leqno]{amsart}
\fi
\usepackage{amsfonts}
\usepackage{amssymb}
\usepackage{amscd}
\usepackage{latexsym}
\usepackage{verbatim}
\ifusepdf
\usepackage{hyperref}
\else\fi
\ifusepsfont
\usepackage[T1]{fontenc}
\usepackage{times}
\else\fi

\iftwelvept
\else\fi

\newtheorem{Theorem}{Theorem}[section]
\newtheorem{Proposition}[Theorem]{Proposition}
\newtheorem{Lemma}[Theorem]{Lemma}
\newtheorem{Corollary}[Theorem]{Corollary}
\newtheorem{Claim}{Claim}[Theorem]

\theoremstyle{definition}

\newtheorem{Remark}[Theorem]{Remark}
\newtheorem{Example}[Theorem]{Example}

\newtheorem{Assertion}[Theorem]{Assertion}

\renewcommand{\theTheorem}{\arabic{section}.\arabic{Theorem}}

\def\rom{\textup}
\newcommand{\ZZ}{{\mathbb{Z}}}
\newcommand{\QQ}{{\mathbb{Q}}}
\newcommand{\RR}{{\mathbb{R}}}
\newcommand{\CC}{{\mathbb{C}}}
\newcommand{\PP}{{\mathbb{P}}}
\newcommand{\OO}{{\mathcal{O}}}
\newcommand{\Aff}{{\mathbb{A}}}

\newcommand{\Spec}{\operatorname{Spec}}
\newcommand{\Pic}{\operatorname{Pic}}
\newcommand{\Cl}{\operatorname{Cl}}
\newcommand{\ord}{\operatorname{ord}}
\newcommand{\Div}{\operatorname{Div}}
\newcommand{\Supp}{\operatorname{Supp}}

\newcommand{\Amp}{\operatorname{Amp}}
\newcommand{\Nef}{\operatorname{Nef}}

\newcommand{\Hilb}{{\operatorname{Hilb}}}
\newcommand{\Hom}{{\operatorname{Hom}}}

\newcommand{\calH}{{\mathcal{H}}}

\newcommand{\id}{\operatorname{id}}
\newcommand{\Image}{\operatorname{Im}}

\newcommand{\Gal}{\operatorname{Gal}}

\newcommand{\GL}{\operatorname{GL}}
\newcommand{\SL}{\operatorname{SL}}

\newcommand{\topo}{\operatorname{top}}

\newcommand{\Proof}{{\sl Proof.}\quad}

\newcommand{\QED}{{\unskip\nobreak\hfil\penalty50\quad\null\nobreak\hfil
{$\Box$}\parfillskip0pt\finalhyphendemerits0\par\medskip}}
\newcommand{\rest}[2]{\left.{#1}\right\vert_{{#2}}}

\begin{document}

\title[surface automorphisms]%
{Projective surface automorphisms of positive                  
topological entropy from an arithmetic viewpoint}
\author{Shu Kawaguchi}
\address{Department of Mathematics, Faculty of Science,
Kyoto University, Kyoto, 606-8502, Japan}
\email{kawaguch@math.kyoto-u.ac.jp}
\begin{abstract}
Let $X$ be a smooth projective surface over a number field 
$K (\subset \CC)$, and $f: X \to X$ an automorphism of positive topological 
entropy. In this paper, we show that there are only 
finitely many $f$-periodic curves on $X$. 
Then we define a height function $\widehat{h}_D$ 
corresponding to a certain nef and big $\RR$-divisor $D$ on $X$ 
and transforming well relative to $f$, and deduce some arithmetic 
properties of $f$-periodic points and non $f$-periodic points.
\end{abstract} 

\maketitle

\section*{Introduction}
\renewcommand{\theTheorem}{\Alph{Theorem}} 
Let $X$ be a smooth projective surface over a number field 
$K \subset \CC$, and $f: X \to X$ an automorphism. 

For $x \in X(\overline{K})$, we denote by  
$O_{f}^+(x) := \{f^n(x) \mid n \geq 1\}$ the forward orbit of $x$ 
under $f$. We fix an ample line bundle $H$ on $X$ and take a 
height function $h_H: X(\overline{K})\to \RR$ corresponding to $H$. 
A point $x$ is said to be {\em $f$-periodic} if $f^n(x) = x$ for some 
positive integer $n$. Since $f$ is an automorphism, $x$ is 
$f$-periodic if and only if $O_{f}^+(x)$ is a finite set. 
If $O_{f}^+(x)$ is not a finite set, we estimate an arithmetic size 
of $O_{f}^+(x)$ by defining   
\[
  N^+_{f,x}(T) := \# \{y \in O_{f}^+(x) \mid h_H(y) \leq T\}
\]
for $T >0$. 

We set $X_{\CC} := X \times_{\Spec(K)} \Spec(\CC)$, and write 
$f_{\CC}: X_{\CC} \to X_{\CC}$ for the induced automorphism. 
Then $f_{\CC}$ induces the linear transformation $f_{\CC}^{*}$ 
on $H^{1,1}(X_{\CC}, \RR)$. Let $\lambda$ be the spectral radius 
of $f_{\CC}^{*}$, i.e., 
the maximum of the absolute values of all the eigenvalues 
of $f_{\CC}^{*}$. Note that $\lambda$ is equal to the 
{first dynamical degree} $\lambda_{1}(f_{\CC})$ of $f_{\CC}$, 
defined by $\lambda_{1}(f_{\CC}) 
:= \lim_{n\to+\infty} \Vert f_{\CC}^{n*} \Vert^{\frac{1}{n}}$ 
(see for example \cite[\S1.4]{DF} and \cite{DS}). 

In this paper, we prove the following assertion (a more precise form 
will be given below), in a sense generalizing \cite{SiK3}. 

Recall that a subset $S \subset X(\overline{K})$ is called 
a set of bounded height if there is a constant $c$ 
such that $h_{H}(s) \leq c$ for all $s \in S$.

\begin{Assertion}
\label{assertion}
Assume $\lambda >1$. Then 
\begin{enumerate}
\item[(1)]
there are a height function $\widehat{h}$ on $X(\overline{K})$ 
and a proper Zariski-closed subset $Z_1$ of $X$ such that 
$\widehat{h}$ transforms well relative to $f$ and such that 
$\widehat{h}$ satisfies the Northcott finiteness property  
on $(X\setminus Z_1)(\overline{K})$. 
\item[(2)]
There is a proper Zariski-closed subset $Z_2$ of $X$ such that 
\[
  \{ x \in (X\setminus Z_2)(\overline{K}) \mid 
  \text{$x$ is $f$-periodic}\}
\]
is a set of bounded height. 
\item[(3)]
There is a proper Zariski-closed subset $Z_3$ of $X$ such that, 
if $x \in \left(X\setminus Z_3 \right)(\overline{K})$ 
is a non $f$-periodic point, then 
\begin{equation}
\label{eqn:assymptotic}
  \lim_{T\to+\infty}\frac{N^+_{f,x}(T)}{\log T} 
  = \frac{1}{\log\lambda}. 
\end{equation}
\end{enumerate}
\end{Assertion}

Before discussing Assertion~\ref{assertion} in more detail, 
we remark that (a similar) assertion is known to be true for some 
other $(X, f)$. 
We consider a modified assertion of 
Assertion~\ref{assertion}, in which $\lambda$ is replaced by the first
dynamical degree $\lambda_{1}(f_{\CC})$, and $f$-periodicity by
$f$-preperiodicity. (Here $x$ is said to be $f$-preperiodic if
$O_f^+(x)$ is a finite set.) 

The first case we remark is when $X$ is a
smooth projective variety and $f: X \to X$ is a morphism such that
$f^*(L) \sim d L$ for some ample divisor $L$ and $d >1$. This case is
studied in Call--Silverman \cite{CS}. Here the modified assertion
holds true with $Z_1 = Z_2 = Z_3 = \emptyset$ by \cite[Theorem~1.1 and
Corollary~1.1]{CS}.  We note that one can check $\lambda_1(f_{\CC}) =
d \;(> 1)$ in this case, using $\lambda_{1}(f_{\CC}) =
\lim_{n\to+\infty}\left( \int_X f^{n*}\omega \wedge \omega^{\dim
    X_{\CC}-1}\right)^{\frac{1}{n}}$ for any normalized K{\"a}hler
form $\omega$ on $X_{\CC}$ (\cite[p.~314]{DS} and 
\cite[Remark~1.17]{DF}).  
The second case is when $X$ is $\PP^2$ and
$f: \PP^2 \dasharrow \PP^2$ is a birational extension of a polynomial
automorphism of the affine plane $f_0: \Aff^2 \to \Aff^2$. Then the
modified assertion holds true with $Z_1 = Z_2 = Z_3 = H_{\infty}$,
where $H_{\infty}$ is the hyperplane at infinity, i.e., $H_{\infty} =
\PP^2\setminus \Aff^2$.  For details, see \cite{Ka2} and the
references therein.  We warn, however, that there are $(X, f)$ for
which the modified assertion does not hold true. (For example,
consider a morphism $f = (f_1, f_2): \PP^1 \times \PP^1 \to \PP^1
\times \PP^1$, where $f_1: \PP^1 \to \PP^1$ is a linear map and $f_2:
\PP^1 \to \PP^1$ is a quadratic map.)

\smallskip
Let us go back to our setting where $f$ is an automorphism of a smooth 
projective surface $X$. Since the works of Gromov and Yomdin state that 
$\log\lambda$ is equal to the topological entropy 
$h_{\topo}(f_{\CC})$ of $f_{\CC}$, the condition $\lambda > 1$ in 
Assertion~\ref{assertion} means that $f_{\CC}$ has 
positive topological entropy. 
It is known that, from a dynamical viewpoint, automorphisms of 
positive topological entropy enjoy many interesting properties  
(cf. \cite{Ca1, Ca, DF, DS, Mc, Og}). 
According to Cantat
\cite[Proposition~1]{Ca1}, if $X$ has an automorphism of positive
topological entropy, then $X$ is either a non-minimal rational
surface, or a surface birational to a $K3$ surface, an Enriques
surface or an abelian surface. 

To single out $Z_1, Z_2$ and $Z_3$ in Assertion~\ref{assertion}, 
we first show the following.  
Here a curve means a reduced 
effective divisor on $X_{\CC}$, and an irreducible curve $C$ is said to be
{\em $f_{\CC}$-periodic} if $f_{\CC}^n(C) = C$ 
for some positive integer $n$. 

\begin{Proposition}[cf. Proposition~\ref{proposition:finiteness}]
\label{prop:A}
Let $f_{\CC}: X_{\CC}\to X_{\CC}$ be an automorphism of 
a smooth projective surface, which has positive topological
entropy. Then there are only finitely many 
$f_{\CC}$-periodic curves on $X_{\CC}$.
\end{Proposition}

In fact, we construct a nef and big $\RR$-divisor class $\nu$ in the
N{\'e}ron-Severi group tensored by $\RR$ such that an irreducible curve
$C$ on $X_{\CC}$ is $f_{\CC}$-periodic if and only if $(\nu, [C]) = 0$.  
By Proposition~\ref{prop:A}, we define the maximal $f_{\CC}$-invariant 
curve $E$ by $E := \sum C$, where $C$ runs all $f_{\CC}$-periodic curves 
on $X_{\CC}$. 
As we will see later, we can take $Z_1 = Z_2 = Z_3 = E$. 
We remark that $E=\emptyset$ if and only if $\nu$ is ample, 
and, in case $E\neq\emptyset$, general points on $E$ behaves differently 
from general points on $X\setminus E$ 
(Theorem~\ref{thm:non-periodic:points}). 
We also remark that, as a corollary of Proposition~\ref{prop:A}, we see 
that an automorphism $f_{\CC}$ of positive topological 
entropy has a Zariski-dense orbit (cf. Corollary~\ref{Cor:og} and 
\cite[Theorem~1.4(2)]{Og}). 

In the following, for simplicity, we say that $f$ has positive topological 
entropy if $f_{\CC}$ does, and we use the notation $h_{\topo}(f)$ instead of 
$h_{\topo}(f_{\CC})$. We show in  Lemma~\ref{lemma:D} that 
we can take an $\RR$-divisor $D$ with $[D] = \nu$ such that 
$f^*(D) + f^{-1*}(D)$ is $\RR$-linearly equivalent to 
$(\lambda + \lambda^{-1}) D$. 
In our setting where $f:X \to X$ is defined over a number field 
$K \subset \CC$, $E$ is defined over $\overline{K}$ 
and $D$ can be taken so that it is defined over $\overline{K}$ 
(Lemma~\ref{lemma:field:def}). 

The following is a precise form of Assertion~\ref{assertion}(1). 

\begin{Theorem}[cf. Theorem~\ref{thm:canonical:height}]
\label{thm:B}
Let $X$ be a smooth projective surface over a number field $K \subset \CC$, 
and $f: X \to X$ an automorphism of positive 
topological entropy $h_{\topo}(f) = \log\lambda$. 
Let $D$ be the nef and big $\RR$-divisor as above, and 
$E$ the maximal $f$-invariant curve. Replacing $K$ by a suitable 
finite extension if necessary, we assume that $D$ and $E$ are defined over 
$K$. Then there exists a unique function 
\[
  \widehat{h}_{D}: X(\overline{K}) \to \RR
\]
with the following two properties. 
\begin{enumerate}
\item[(i)]
$\widehat{h}_D$ is a height function corresponding 
to $D$.  
\item[(ii)]
$\widehat{h}_D(f(x)) + \widehat{h}_D(f^{-1}(x))
= \left(\lambda + \lambda^{-1}\right) \widehat{h}_D(x)$ 
for all $x \in X(\overline{K})$.  
\end{enumerate} 
Moreover, $\widehat{h}_{D}$ satisfies the following properties. 
\begin{enumerate}
\item[(iii)] 
$\widehat{h}_D$ is non-negative. 
\item[(iv)]
For $x \in X(\overline{K})$, $\widehat{h}_D(x) = 0$ if and only 
if $x$ lies on $E$ or $x$ is $f$-periodic.  
\item[(v)]
$\widehat{h}_D$ satisfies the Northcott finiteness property 
on $(X\setminus E)(\overline{K})$ 
\textup{(}See Theorem~\ref{thm:canonical:height} 
for a precise statement\textup{)}. 
\end{enumerate} 
We call $\widehat{h}_{D}$ a 
{\em canonical height function} for $f$ and $D$. 
\end{Theorem}

In a pioneering paper \cite{SiK3}, Silverman constructed canonical height
functions $\widehat{h}$ on $K3$ surfaces given by the complete
intersection of $(1,1)$ and $(2,2)$ hypersurfaces in $\PP^2 \times
\PP^2$, and showed that $\widehat{h}$ enjoys properties comparable 
with the N{\'e}ron-Tate height functions on abelian varieties.   
These $K3$ surfaces have two involutions $\sigma_1, \sigma_2$ induced from  
the projections $p_1, p_2: X \to \PP^2$. 
One sees that the composition $\sigma_2\circ\sigma_1$
has the positive topological entropy $\log(7+4\sqrt{3})$, and
$\widehat{h}$ is the canonical height function for 
$\sigma_2\circ\sigma_1$ and $(1 +\sqrt{3})(p_1^* H + p_2^* H)$, 
where $H$ is a line in $\PP^2$.   
From this view,
Theorem~\ref{thm:B} gives an extension of \cite[Theorem~1.1]{SiK3} 
to any smooth projective surface automorphisms of positive topological 
entropy. We also note that Theorem~\ref{thm:B}, together with 
Theorem~\ref{thm:C} below, gives a refinement of \cite[\S1]{Ka} 
when $(X; f, f^{-1})$ is viewed as a dynamical eigensystem 
for $D$ of degree $(\lambda + \lambda^{-1})$.

On the above $K3$ surfaces, there are no $(\sigma_2\circ\sigma_1)$-periodic 
curves and so $D$ is ample. In general, 
projective surface automorphisms of positive 
topological entropy have non-empty $E$ and so $D$ cannot be ample. 
Still, there is an effective divisor $Z$ on $X$ such that 
$\Supp(Z) \subseteq \Supp(E)$ and such that $D - \varepsilon Z$ 
is ample for sufficiently small $\varepsilon > 0$.  
Moreover, since $E$ is a curve, the dynamical properties of $f$ restricted 
to $E$ are easy to understand. 
By these, together with Theorem~\ref{thm:B}, 
we get Assertion~\ref{assertion}(2)(3). 

\begin{Theorem}[cf. Theorem~\ref{thm:periodic:points}, 
Corollary~\ref{cor:non-periodic:points}]
\label{thm:C}
With the notation and assumption in Theorem~\ref{thm:B}, 
we have the following. 
\begin{enumerate}
\item[(1)] 
The set $\{x \in (X\setminus E)(\overline{K})\mid 
\rom{$x$ is $f$-periodic}\}$ is of bounded height. 
\item[(2)] 
Let $x \in (X\setminus E)(\overline{K})$ be a non $f$-periodic point. 
Then 
\[
  \lim_{T\to+\infty}\frac{N^+_{f,x}(T)}{\log T} 
  = \frac{1}{\log\lambda}.
\] 
\textup{(}Also see Theorem~\ref{thm:non-periodic:points} 
for stronger statements.\textup{)}
\end{enumerate}
\end{Theorem}

The organization of this paper is as follows. 
After reviewing some basic facts on smooth projective 
surfaces over $\CC$ in \S\ref{sec:preliminary}, 
we construct a nef and big $\RR$-divisor class $\nu$ in \S\ref{sec:nef}, 
with which we derive finiteness of $f$-periodic curves 
in \S\ref{sec:periodic}. Also in \S\ref{sec:periodic}, 
we give examples of automorphisms of positive topological entropy, 
and construct a nef and big $\RR$-divisor $D$ with $[D] = \nu$.  
Now turning our attention to arithmetic properties, 
we briefly review on height functions in \S\ref{sec:height}. 
In \S\ref{sec:canonical}, we construct canonical height functions 
$\widehat{h}_{D}$ and study their further properties.  
In \S\ref{sec:thm:B:C}, we prove Theorem~\ref{thm:C} and related results. 

One approach to prove Theorems~\ref{thm:B} and \ref{thm:C} would be 
to contract the curve $E$ on $X$ and to work over 
the singular surface $Y$ on which $E$ is contracted. In general, however, 
this is possible as analytic surfaces (cf. \cite[Example~5.7.3]{Ha}), 
and it is not clear to the author whether in our case $Y$ is algebraic 
(which is preferably $\QQ$-factorial). 
So, we choose to work over $X$, using   
the ampleness of $D - \varepsilon Z$. 

Cantat \cite{Ca1} showed the existence of the unique invariant
probability measure $\mu$ of maximal entropy on $X$. 
Finally let us note that it would be an interesting question 
whether or not the Galois orbits of general periodic points 
(or of general points of small height) are
asymptotically equidistributed with respect to $\mu$ 
(cf. \cite[Question~3.4.1]{Ka}).

\smallskip
{\it Acknowledgments.}\quad
I thank Professors Takeshi Abe, Vincent Maillot,   
Shigeru Mukai and Keiji Oguiso for valuable discussions and assistance. 
I thank the referee for carefully reading the manuscript and
giving many helpful comments. 

\medskip
\renewcommand{\theTheorem}{\arabic{section}.\arabic{Theorem}}
\setcounter{equation}{0}
\section{Preliminaries on projective surfaces}
\label{sec:preliminary}
In this section, we review some facts on projective 
surfaces and fix notation that will be used in the sequel. 
We refer the reader to \cite{BHPV}, \cite{Bea} and \cite{Laz} for details.   

Let $X$ be a smooth projective surface over $\CC$. 
A curve means a reduced effective divisor on $X$. 
We write for $\Div(X)$  
the group of divisors on $X$, and 
$\Cl(X)$ for the divisor class group on $X$, i.e., 
the linearly equivalence classes of divisors on $X$. 
Let $N^1(X)$ be the N{\'e}ron-Severi group, i.e., the numerically
equivalence classes of divisors on $X$. We denote its rank by $\rho$. 
We write $\Div_{\RR}(X) := \Div(X) \otimes\RR$ for the $\RR$-vector space 
of $\RR$-divisors on $X$, and $\Cl(X)_{\RR}$ for the $\RR$-vector space 
of $\RR$-linearly equivalence classes of $\RR$-divisors on $X$, 
and $N^1(X)_{\RR}$ for the $\RR$-vector space 
of numerically equivalence classes of $\RR$-divisors on $X$. 
We note that $\Cl(X)_{\RR} = \Cl(X)\otimes\RR$ and 
$N^1(X)_{\RR} = N^1(X) \otimes \RR$.  
For $\RR$-divisors $D_1, D_2$, we write 
$D_1 \sim_{\RR} D_2$ if $D_1, D_2$ are $\RR$-linearly 
equivalent. For an $\RR$-divisor $D$, 
we denote by $[D]$ its image in $N^1(X)_{\RR}$. 
An $\RR$-divisor $D$ is said to be integral if it lies in 
the image of the natural map $\Div(X) \to \Div(X)_{\RR}$. 

An $\RR$-divisor $D$ is said to be {\em ample} (resp.\ {\em big}, 
{\em effective}) if $D$ is expressed as 
a finite sum $D = \sum c_i D_i$ for some 
$c_i > 0$ and ample (resp.\ big, effective) integral divisors $D_i$. 
An $\RR$-divisor $D$ is said to be {\em nef} if $(D, C) \geq 0$ 
for all irreducible curves $C$ on $X$. 
Amplitude, bigness and nefness depend only on 
numerically equivalence classes. 
Let $\Amp(X) \subset N^1(X)_{\RR}$ and 
$\Nef(X) \subset N^1(X)_{\RR}$ denote the ample cone 
and the nef cone in $N^1(X)_{\RR}$, respectively. 
Kleiman's theorem states that the closure of the ample 
cone is equal to the nef cone (\cite[Theorem~1.4.23]{Laz}): 
\begin{equation}
\label{eqn:Kleiman}
\overline{\Amp(X)} = \Nef(X). 
\end{equation}

We collect some facts on $\RR$-divisors. For a proof, see for example  
\cite[Proposition~2.2.22, the proof of Theorem~2.3.18]{Laz}. 

\begin{Proposition}
\label{prop:divisor}
Let $D$ be an $\RR$-divisor on $X$. 
\begin{enumerate}
\item[(1)]
$D$ is big if and only if $[D] = [A] + [N]$ in $N^1(X)_{\RR}$ , where 
$A$ is an ample $\RR$-divisor and $N$ is an effective $\RR$-divisor. 
\item[(2)]
Assume $D$ is nef. Then 
$D$ is big if and only if $(D^2) >0$. 
\end{enumerate}
\end{Proposition}

Since $N^1(X)_{\RR}$ has signature $(1, \rho-1)$ by the Hodge index 
theorem, the set $\{x \in N^1(X)_{\RR} \mid (x,x) >0 \}$ consists 
of two disjoint cones. Let $P$ be the positive cone, i.e., 
one of the cones that contains an ample divisor. 
Let $\overline{P}$ be the closure of $P$ in $N^1(X)_{\RR}$. 
The following elementary lemma will be useful. 
See for example \cite[IV,~Cor.~7.2]{BHPV} about (1) and (2). 
The Cauchy--Schwartz inequality yields (3).

\begin{Lemma}
\label{lemma:useful}
\begin{enumerate}
\item[(1)] If $x, y \in \overline{P}$, then $(x,y)\geq 0$. 
\item[(2)] If $x \in \overline{P}\setminus\{0\}$ and $y \in P$, 
then $(x, y) > 0$. 
\item[(3)] If $x, y \in \overline{P}\setminus\{0\}$ and 
$x \neq t y$ for all $t \in \RR_{>0}$, then $(x, y) > 0$. 
\end{enumerate}
\end{Lemma}

The next proposition will be used in \S\ref{sec:periodic} and 
\S\ref{sec:canonical}. 

\begin{Proposition}
\label{prop:nu:2}
Let $\nu \in N^1(X)_{\RR}$ be a nef and big $\RR$-divisor class 
on $X$. We set 
\[
  Z(\nu) :=\{C \mid \text{$C$ is an irreducible 
  curve on $X$ such that $(\nu, [C]) = 0$}\}.
\] 
Then 
\begin{enumerate}
\item[(1)]
$Z(\nu)$ is a finite set. 
\item[(2)]
There is an effective divisor $Z$ on $X$ such that 
$\Supp(Z) \subseteq \bigcup_{C \in Z(\nu)} \Supp(C)$ 
and such that $\nu - \varepsilon [Z]$ is ample for all sufficiently 
small $\varepsilon >0$.  
\end{enumerate}
\end{Proposition}

\Proof
(1) From Proposition~\ref{prop:divisor}, there are an ample 
$\RR$-divisor $A$ and an effective $\RR$-divisor $N$ such 
that $\nu = [A] +[N]$. If an irreducible curve $C$ satisfies 
$(\nu, [C]) =0$, then $(N, C) = - (A, C) < 0$. Hence $C$ must 
be an irreducible component of $\Supp(N)$. Since the number of 
irreducible components of $\Supp(N)$ is finite, $Z(\nu)$ is a finite 
set.
 
\smallskip
(2) 
We set $Z(\nu) = \{C_1, \ldots, C_m\}$.  
Since $(\nu^2) > 0$, 
the Hodge index theorem asserts that 
the intersection matrix $((C_i, C_j))_{i,j} \leq 0$.  
Then there is an effective divisor 
$Z_1 = \sum_{i=1}^m a_i C_i$ such that, for all $i$, 
$(Z_1, C_i) \leq 0$ and $a_i >0$ (cf. \cite[Chap.~I~(2.10)]{BHPV}). 

We claim that $\nu - \varepsilon [Z_1]$ is nef and big 
for all sufficiently small $\varepsilon >0$. Indeed, since 
$(\nu^2) >0$, we have $((\nu - \varepsilon [Z_1])^2) > 0$ 
for $1 \gg \varepsilon >0$. 
Let $C$ be an irreducible curve on $X$. 
If $C \in Z(\nu)$, then 
$(\nu - \varepsilon [Z_1], [C]) = - \varepsilon (Z_1, C) \geq 0$. 
If $C \subseteq \Supp(N)$ but $C \not\in Z(\nu)$, then  
$(\nu, [C]) > 0$, so that $(\nu - \varepsilon [Z_1], [C]) > 0$ 
for $1 \gg \varepsilon >0$, 
On the other hand, since $Z(\nu) \subseteq \Supp(N)$, 
$N - \varepsilon Z_1$ is effective for $1 \gg \varepsilon >0$ 
and $\Supp(N - \varepsilon Z_1) \subseteq \Supp(N)$. 
Thus, if $C \not\subseteq \Supp(N)$, then 
$(\nu - \varepsilon [Z_1], [C]) > (N - \varepsilon Z_1, C) \geq 0$. 
We conclude that $\nu - \varepsilon [Z_1]$ is nef and big 
for $1 \gg \varepsilon >0$. 
There is some $C_i$ with $(Z_1, C_i) < 0$. 
(Indeed, if $(Z_1, C_i) = 0$ for all $i$, then $Z_1 = \sum_{i=1}^m a_i 
C_i$ is numerically trivial by the Hodge index theorem, which is absurd 
because $a_i >0$ for all $i$.) Hence  
\[
  Z(\nu - \varepsilon [Z_1]) 
  := \{C \mid \text{$C$ is an irreducible 
  curve on $X$ such that $(\nu - \varepsilon [Z_1], [C]) = 0$}\}
\] 
is a proper subset of 
$Z(\nu)$. Inductively, we apply the above argument 
to $\nu - \varepsilon [Z_1]$, we see that there is 
an effective divisor $Z$ on $X$ with  
$\Supp(Z) \subseteq \bigcup_{C \in Z(\nu)} \Supp(C)$ such that 
$\nu - \varepsilon [Z]$ is nef and big for $1 \gg \varepsilon >0$ 
and such that 
\[
  Z(\nu - \varepsilon [Z])
  := 
  \{C \mid \text{$C$ is an irreducible 
  curve on $X$ such that $(\nu - \varepsilon [Z], [C]) = 0$}\}
\] 
is empty. Then the Nakai criterion for $\RR$-divisors 
due to Campana and Peternell (\cite[Theorem~2.3.18]{Laz})
yields that $\nu - \varepsilon [Z]$ is ample 
for $1 \gg \varepsilon >0$. 
\QED

\medskip
\setcounter{equation}{0}
\section{Nef and big $\RR$-divisor classes}
\label{sec:nef}
Let $X$ be a smooth projective surface over $\CC$, 
and $f: X \to X$ an automorphism. 
Then $f$ induces $f^*: H^2(X, \ZZ) \to H^2(X, \ZZ)$. 
By a slight abuse of notation, we also denote by $f^*$ the induced 
transformations on $H^2(X, \RR)$ and on $H^2(X, \CC)$. 
Let 
\[
  \lambda := \sup 
  \{|t| \mid \det(t I - f^*) = 0\}
\]
be the spectral radius of $f^*$. 

By the works of Gromov and Yomdin, the topological entropy of $f$ 
is related to the spectral radius of $f^*$, which we might take 
as one of the definitions of topological entropy. 
We put together in Theorem~\ref{thm:lambda} the 
works of Gromov and Yomdin and 
other results which will be needed in the sequel. 
For Theorem~\ref{thm:lambda}(1), 
we refer the reader to 
\cite[Theorem~2.1.5]{Ca} and \cite[Th{\'e}or{\`e}me~2.1]{DS}.   
For Theorem~\ref{thm:lambda}(2), 
we refer the reader to \cite[Theorem~3.2]{Mc}. We note that 
$K3$ surfaces are treated in \cite[Theorem~3.2]{Mc}, 
but since the signature of 
$H^{1,1}(X, \RR)$ is $(1, \rho-1)$, the argument there also applies to 
automorphisms of smooth projective surfaces. 

\begin{Theorem}[cf. \cite{Ca, DS, Mc}]
\label{thm:lambda}
Let $f: X \to X$ be an automorphism of a 
smooth projective surface on $X$. 
\begin{enumerate}
\item[(1)]
The topological entropy $h_{\topo}(f)$ is equal to $\log\lambda$. 
Moreover, $h_{\topo}(f)$ is also equal to the spectral radius of 
$\rest{f^*}{H^{1,1}(X, \RR)}$ acting on $H^{1,1}(X, \RR)$. 
\item[(2)]
Assume $f$ has positive topological entropy, i.e., 
$\lambda > 1$. Then the set of eigenvalues of $f^*$ 
with counted multiplicities is given by 
\[
  \{\lambda, \lambda^{-1}, 
  \alpha_1, \alpha_2, \cdots, \alpha_{\dim H^2(X, \RR)-2}\},
\] 
where $|\alpha_i| = 1$ for all $i = 1, 2, \ldots, \dim H^2(X, \RR) -2$. 
In particular, $f^*$ has a unique, simple eigenvalue $\beta$
with $|\beta| >1$, and $\beta$ is equal to $\lambda$. 
\end{enumerate}
\end{Theorem}

\begin{Lemma}
\label{lemma:contraction}
Let $X$ be a smooth projective surface over $\CC$, 
and $f: X \to X$ an automorphism. We assume that 
the Kodaira dimension of $X$ is non-negative.  
Let $X_0$ be the minimal surface of $X$, and 
$\pi: X \to X_0$ the contraction morphism. 
Then there is a unique automorphism 
$f_0: X_0 \to X_0$ with  
$\pi\circ f = f_0 \circ \pi$. Moreover, 
$h_{\topo}(f_0)= h_{\topo}(f)$. 
\end{Lemma}

\Proof
The automorphism $f: X \to X$ induces a birational map 
$f_0: X_0 \dasharrow X_0$ with $\pi\circ f = f_0 \circ \pi$. 
Since the Kodaira dimension of $X$ is non-negative, 
$f_0$ is an automorphism (\cite[Theorem~V.19]{Bea}). 

Let $\{F_1, \ldots, F_k\}$ be the set of exceptional irreducible 
curves of $\pi$, i.e., the set of irreducible curves $F_i$'s that 
are contracted to points by $\pi$. Then 
$H^2(X, \ZZ) = \pi^*(H^2(X_0, \ZZ)) 
\oplus \bigoplus_{i=1}^k \ZZ [F_i]$. 

Since $\pi\circ f = f_0 \circ\pi$, 
$f^*$ preserves $\pi^*(H^2(X_0, \ZZ))$ and 
$\bigoplus_{i=1}^k \ZZ [F_i]$. 
On the other hand, since the intersection form on 
$\bigoplus_{i=1}^k \ZZ [F_i]$ is negative definite, 
Lemma~\ref{lemma:Mc} below asserts that 
the eigenvalues of $\rest{f^*}{\bigoplus_{i=1}^k \ZZ [F_i]}$ 
lie on the unit circle. Hence 
$f_0$ has the same topological entropy as $f$ 
by Theorem~\ref{thm:lambda}(1). 
\QED

\begin{Lemma}[\cite{Mc}, Lemma~3.1]
\label{lemma:Mc}
Let $O(p,q) \subset \GL_{p+q}(\RR)$ be the orthogonal group 
of the real quadratic form 
$x_1^2 + \cdots + x_p^2 - x_{p+1}^2 - \cdots - x_{p+q}^2$ 
of signature $(p,q)$. Then, for any $T \in O(p,q)$, 
$T$ has at most $\min(p,q)$ eigenvalues outside the unit circle, 
counted with multiplicities.  
\end{Lemma}

Cantat \cite[Proposition~1]{Ca1} classified 
smooth projective surfaces that have automorphisms of 
positive topological entropy (cf. Lemma~\ref{lemma:contraction}).  

\begin{Proposition}[\cite{Ca1}, Proposition~1]
\label{prop:Ca}
Let $X$ be a smooth projective surface over $\CC$ which has an 
automorphism of positive topological entropy. Then 
$X$ is either a non-minimal rational surface, or 
a surface birational to a $K3$ surface, an Enriques surface or 
an abelian surface. 
\end{Proposition}

From now on, we consider automorphisms of 
positive topological entropy. 
The next proposition concerns the existence of a nef 
and big $\RR$-divisor class $\nu$ that transforms well 
relative to $f$. This may be seen as an $\RR$-divisor class version of 
\cite[Th{\'e}or{\`e}me~2]{Ca1} and \cite[Theorem~5.1]{DF}. 

\begin{Proposition}
\label{prop:nu}
Let $X$ be a smooth projective surface over $\CC$, 
and $f: X \to X$ an automorphism of positive 
topological entropy $\log\lambda$. 
\begin{enumerate}
\item[(1)]
There are non-zero nef classes $\nu_+$ and $\nu_- \in \Nef(X)$ 
such that $f^*(\nu_+) = \lambda \nu_+$ and 
$f^*(\nu_-) = \lambda^{-1} \nu_-$. 
\item[(2)]
We have $(\nu_+^2) = 0$ and $(\nu_-^2) = 0$. 
If $\nu^{\prime}_+, \nu^{\prime}_- \in \Nef(X)$  
are other nef classes with $f^*(\nu^{\prime}_+) = \lambda \nu^{\prime}_+$, 
$f^*(\nu^{\prime}_-) = \lambda^{-1} \nu^{\prime}_-$, 
then $\nu^{\prime}_+ = c_+ \nu_+$, $\nu^{\prime}_- = c_- \nu_-$
for some $c_+, c_- > 0$. 
\item[(3)]
Set $\nu := \nu_+ + \nu_- \in N^1(X)_{\RR}$. 
Then $\nu$ is nef and big. Moreover, one has 
\[
  f^*(\nu) + f^{-1*}(\nu) = 
  (\lambda + \lambda^{-1}) \nu. 
\]
\end{enumerate}
\end{Proposition}

\Proof
(1)(2) 
This is remarked in \cite[Remerque~1.1]{Ca} for $K3$ surfaces, 
whose argument also applies to other surfaces. 
For the sake of completeness, we briefly give a proof. 

By Theorem~\ref{thm:lambda}(2), there are non-zero $\nu_+$ and  
$\nu_- \in H^2(X, \RR)$ such that $f^*(\nu_+) = \lambda \nu_+$  
and $f^*(\nu_-) = \lambda^{-1} \nu_-$. 
We claim that 
the classes $\nu_+$ and $\nu_-$ belong to $N^1(X)_{\RR}$. 
Indeed, we set 
\[
  W := \{x \in H^2(X, \RR) \mid 
  \text{$(x, \delta) = 0$ for all $\delta \in N^1(X)$}\}.
\] 
If $X$ is a rational surface or a surface birational to 
an Enriques surface, then $W=0$, so that we get the claim trivially. 
If $X$ is a surface birational to an $K3$ surface or 
an abelian surface, then the Lefschetz $(1,1)$-theorem 
tells us that $\rest{f^*}{W}$ is conjugate to an element of 
$O(2, 0) \times O(0, \dim H^2(X, \RR) - 2 - \rho)$. 
Then, from Lemma~\ref{lemma:Mc}, all eigenvalues of $\rest{f^*}{W}$ 
lie on the unit circle. Thus we get the claim in this case also.  

Now Proposition~\ref{prop:nu}(1)(2) is proven as 
in \cite[Th{\'e}or{\`e}me~2]{Ca1} and \cite[Theorem~5.1]{DF}. 
Noting that $\Amp(X)$ spans $N^1(X)_{\RR}$ as an $\RR$-vector space, 
we take a general element $\gamma$ in $\Amp(X)$ and fix it. 
Since $\lambda$ is the unique, simple eigenvalue of $f^*$ 
whose absolute value is greater than $1$ by 
Theorem~\ref{thm:lambda}(2), 
$\frac{1}{\lambda^n} (f^*)^n(\gamma)$ converges to $c \nu_+$ 
with some $c \neq 0$ as $n \to +\infty$. 
We replace $\nu_+$ by $c \nu_+$. 
Since $\frac{1}{\lambda^n} (f^*)^n(\gamma) \in \Amp(X)$ for all 
$n \geq 1$, we see that $\nu_+ \in \Nef(X)$ by Kleiman's theorem 
\eqref{eqn:Kleiman}. A similar argument shows that 
$\nu_-$ can also be taken in $\Nef(X)$. 

It follows from $\lambda^2 (\nu_+^2) = (f^*(\nu_+)^2) 
= (\nu_+^2)$ that $(\nu_+^2) = 0$. 
Similarly, $(\nu_-^2) = 0$. 
By Theorem~\ref{thm:lambda}(2), $\nu_+$ is unique up to 
positive real numbers and so is $\nu_-$.  

(3) 
Obviously, $\nu = \nu_+ + \nu_-$ is nef.  
Since $(\nu^2) = (\nu_+^2) + 2 (\nu_+, \nu_-) + (\nu_-^2)= 
2 (\nu_+, \nu_-)$, 
Lemma~\ref{lemma:useful}(3) yields that $(\nu^2) > 0$. 
By Proposition~\ref{prop:divisor}(2), 
$\nu$ is big. 
We have
\begin{align*}
  f^*(\nu) + f^{-1*}(\nu) &= 
  f^*(\nu_+) + f^*(\nu_-) + f^{-1*}(\nu_+) + f^{-1*}(\nu_-) \\ 
  & = (\lambda + \lambda^{-1}) (\nu_+ + \nu _-) 
  = (\lambda + \lambda^{-1})\nu. 
\end{align*}
\QED

\begin{Remark}
\label{remark:quadratic}
If $\lambda$ is a quadratic irrational number, then 
$\nu$ can be taken as an integral divisor class, i.e., 
$\nu \in N^1(X)$. Indeed, in this case, $\lambda$ is a root of 
the quadratic polynomial $t^2 - \tau t +1$, where $\tau := \lambda +
\lambda^{-1}$ is an integer.  
We set $F = \QQ(\lambda)$ and let $O_F$
be the ring of integers of $F$. Let $\sigma \in \Gal(F/\QQ)$ be 
the non-trivial element. 

Then there is a nef class $\nu_+ \in N^1(X) \otimes_{\ZZ} O_F$ 
with $f^* \nu_+ = \lambda \nu_+$. We set $\nu_- := \sigma(\nu_+) 
\in N^1(X) \otimes_{\ZZ} O_F$. 
Then $\nu_-$ is nef and $f^*(\nu_-) = \lambda^{-1} \nu_-$. 
We set $\nu := \nu_+ + \nu_-$. Since $\sigma(\nu) = \nu$, we see 
that $\nu \in N^1(X)$. 
\end{Remark}

\medskip
\setcounter{equation}{0}
\section{Finiteness for $f$-periodic curves}
\label{sec:periodic}
Let $X$ be a smooth projective surface over $\CC$, 
and $f: X \to X$ an automorphism of positive topological entropy.  
A point $x$ on $X$ is said to be {\em $f$-periodic} if 
$f^n(x) = x$ for some positive integer $n$. 
An irreducible curve $C$ on $X$ is said to be {\em $f$-periodic} 
if $f^n(C) = C$ for some positive integer $n$. Note that an 
$f$-periodic $C$ may {\em not} be fixed pointwisely by $f$. 

By \cite[Th{\'e}or{\`e}me~in~p.~906]{Ca1}, 
if $f$ has positive topological entropy, 
then there are infinitely many $f$-periodic points on $X$. 
Indeed, suppose there are only finitely many periodic points, 
then the support of the measure of $\mu$ in {\it op.cit.} 
is a finite set. This is impossible since the measure entropy 
with respect to $\mu$ is $\log\lambda >0$.   

In contrast, the next proposition shows that there are only finitely many 
$f$-periodic curves on $X$. 

\begin{Proposition}
\label{proposition:finiteness}
Let $X$ be a smooth projective surface over $\CC$, 
and $f: X \to X$ an automorphism of positive topological entropy. 
\begin{enumerate}
\item[(1)]
There are only finitely many $f$-periodic curves on $X$.
If $X$ is non-rational and $C$ is an $f$-periodic curve, 
then $C$ is isomorphic to $\PP^1$. 
\item[(2)]
Let $\nu \in N^1(X)_{\RR}$ be a nef and big class 
in Proposition~\ref{prop:nu}. 
Then an irreducible curve $C$ is $f$-periodic if and only 
if $([C], \nu) = 0$. Moreover, if 
$X$ is non-rational and $C$ is an irreducible curve with 
$([C], \nu) = 0$, then $C$ is  isomorphic to $\PP^1$. 
\end{enumerate}
\end{Proposition}

\Proof
In virtue of Proposition~\ref{prop:nu:2}, we have only to show (2). 
Let $\log\lambda$ denote the topological entropy of $f$. 

Suppose $C$ is an $f$-periodic curve, i.e., $f^n(C) = C$ for some 
positive integer $n$. For $k \in \ZZ$, we then have 
{\allowdisplaybreaks
\begin{align}
\label{eqn:finiteness}
  ([C], \nu_+) + ([C], \nu_-) 
  & = ([C], \nu) \\ \notag
  & = ([f^{kn}(C)], \nu) = ([C], (f^*)^{kn}(\nu)) \\ \notag
  & = ([C], (f^*)^{kn}(\nu_+)) + ([C], (f^*)^{kn}(\nu_-)) \\ \notag
  & = \lambda^{kn}([C], \nu_+) + \lambda^{-kn}([C], \nu_-). 
\end{align}
}
Since Eqn.\!~\eqref{eqn:finiteness} holds for all $k \in \ZZ$, \
we obtain $([C], \nu_+) = 0$ and $([C], \nu_-) = 0$, 
whence $([C], \nu) = 0$. 

Next, suppose $C$ is an irreducible curve with 
$([C], \nu) =0$. Since $\nu_+$ and $\nu_-$ are nef classes, 
we get $([C], \nu_+) = 0$ and $([C], \nu_-) = 0$. 
Then 
\[
  ([f^k(C)], \nu) = \lambda^{k}([C], \nu_+) + \lambda^{-k}([C], \nu_-) 
  = 0 
\]
for all $k \in \ZZ$. Proposition~\ref{prop:nu:2}(2) yields that 
$f^k(C) = f^l(C)$ for some positive integers $k <l$. 
Thus $f^{l-k}(C) = C$. 

Finally, suppose that $X$ is non-rational. By 
Proposition~\ref{prop:Ca}, $X$ is birational to 
a $K3$ surface, an Enriques surface, or an abelian surface. 
Let $C$ be an $f$-periodic curve. 

{\bf Case 1:} The case $X$ is minimal. 
Since $(\nu^2) >0$ and $(\nu, [C]) =0$, the Hodge index theorem 
yields that $(C^2) <0$. Since $K_{X}$ is numerically trivial, 
we have $(C, K_X)=0$. Then the adjunction formula shows that 
$C$ is a smooth rational curve with selfintersection $-2$. 

{\bf Case 2:} The case $X$ is not minimal. 
Noting that $X$ is non-rational, 
let $X_0$ be the minimal surface of $X$, and $\pi: X \to X_0$ 
the contraction morphism. By Lemma~\ref{lemma:contraction}, 
$f$ induces the automorphism $f_0: X_0 \to X_0$ with 
$h_{\topo}(f_0) = \log\lambda$.  

If $C$ is contracted to a point by $\pi$, 
then $C$ is isomorphic to $\PP^1$. 
If $C$ is not contracted to a point by $\pi$, then 
we set $C_0 = \pi(C)$. Since $\pi\circ f = f_0 \circ \pi$, 
$C_0$ is $f_0$-periodic. Then, by Case 1, 
$C_0$ is isomorphic to $\PP^1$. 
Since $\pi$ is a succession of blow-ups, 
$\rest{\pi}{C}: C\to C_0$ is birational. Hence
$C$ is isomorphic to $\PP^1$. 
\QED

A curve $D$ on $X$ is said to be {\em $f$-invariant}  if $f(D) = D$. 
We define the curve $E$ by 
\begin{equation}
\label{eqn:maximal:curve}
  E = \sum C, 
\end{equation}
where $C$ runs all $f$-periodic curves. 
Then $E$ is $f$-invariant. 
We will call $E$ the {\em maximal $f$-invariant curve}.  
This terminology may be justified from 
the following proposition. 

\begin{Proposition}
\label{prop:E}
The curve $E$ is maximal among $f$-invariant curves. 
Namely, if $D$ is an $f$-invariant curve, 
then $E - D$ is effective.  
\end{Proposition}

\Proof
Let $D = \sum_{i=1}^k D_i$ be the irreducible decomposition 
of $D$. Then $f$ acts on the set $\{D_1, \ldots, D_k\}$ as 
a permutation. In particular, $f^{k!}(D_i) = D_i$ for all $i$. 
This shows that $D_i$'s are $f$-periodic, and thus 
$E - D$ is effective.  
\QED

The following corollary gives a different proof 
of \cite[Theorem~1.4(2)]{Og} (cf. Remark~\ref{rmk:null}), 
when $X$ is an algebraic $K3$ surface.   

\begin{Corollary}
\label{Cor:og}
Let $X$ be a smooth projective surface over $\CC$, 
and $f: X \to X$ an automorphism of positive topological 
entropy. If $x$ is not $f$-periodic and does not lie on $E$, then 
the orbit $\{f^n(x) \mid n \in \ZZ\}$ is Zariski dense in $X$. 
\end{Corollary}

\Proof
Let $W$ be the Zariski closure of $\{f^n(x) \mid n \in \ZZ\}$. 
To derive a contradiction, we assume $W \neq X$. 
Since the one dimensional component $W_1$ of $W$ satisfies 
$f(W_1) = W_1$, we have $W_1 \subseteq \Supp(E)$ by 
Proposition~\ref{prop:E}. Since 
$x$ does not lie on $E$, we conclude $W_1 = \emptyset$. 
Then $W$ becomes a finite set, which contradicts with the assumption 
that $x$ is not $f$-periodic. 
\QED

\begin{Remark}
\label{rmk:null}
In contrast, suppose $X$ is a $K3$ surface and $f: X \to X$ is 
an automorphism of null topological entropy. 
Then $X = \bigcup C$ as a set, where $C$ runs all $f$-periodic 
curves. Indeed, by \cite[Proposition~1.4,~Corollaire~2.2]{Ca}, 
either $f^{n} = \id$ for some positive integer $n$, or 
there is a surjective morphism $\pi: X \to \PP^1$ with $\pi \circ f = \pi$. 
In the first case, every curve is $f$-periodic. In the second case, 
every curve in a fiber is $f$-periodic. 
\end{Remark}

For illustration, we give some examples of 
projective surface automorphisms of positive 
topological entropy. For other examples,  
we refer the reader to \cite{Ca} and \cite{Mc}. 
See also \cite[\S1.4]{Ka} and the references therein. 

\begin{Example}[$K3$ surfaces in $\PP^2 \times \PP^2$]
\label{subsec:eg:1}
Let $X$ be a $K3$ surface in $\PP^2 \times \PP^2$ given by the
intersection of hypersurfaces of bidegree $(1,1)$ and $(2,2)$.  
This surface is sometimes called a Wehler K3 surface (cf. \cite{Ma}).  
When $X$ is defined over a number field $K$, 
Silverman \cite{SiK3} constructed a canonical height function 
on $X(\overline{K})$ and deduced arithmetic properties of $X$. 

For $i = 1, 2$, 
we write $p_i : X \to \PP^2$ for the projection to the $i$-th factor, 
and $\sigma_i: X \to X$ the involution given by $p_i$. 
We take the composition of these two involutions: 
\[
  f = \sigma_2\circ\sigma_1. 
\]

Let $H$ be a hyperplane in $\PP^2$, and 
we set $H_i = p_i^{-1}(H)$. 
As in \cite[\S1]{SiK3}, we put 
\[
  \nu_+ = \left(2+\sqrt{3}\right) [H_1] - [H_2], \quad 
  \nu_+ = -[H_1] + \left(2+\sqrt{3}\right)[H_2] 
\]
in $N^1(X)_{\RR}$. By \cite[Lemma~2.1]{SiK3}, $\nu_+, \nu_-$ 
are nef and satisfy 
\[
  f^*(\nu_+) = \left(7 + 4 \sqrt{3}\right) \nu_+, 
  \quad
  f^*(\nu_-) 
  = \left(7 + 4 \sqrt{3}\right)^{-1} \nu_-. 
\]
The morphism $f$ has the positive topological entropy 
$\log\left(7 + 4 \sqrt{3}\right)$. 
Here $\nu := \nu_+ + \nu_- 
= (1+\sqrt{3}) ([H_1] + [H_2])$ is ample, 
and there are no $f$-periodic curves on $X$ 
(cf. Proposition~\ref{proposition:finiteness}(2)). 

Another example of a $K3$ surface is a smooth hypersurface 
of tridegree $(2,2,2)$ in $\PP^1 \times \PP^1 \times \PP^1$. 
It has three natural involutions $\sigma_1, \sigma_2, 
\sigma_3$, and the composition $f := \sigma_3\circ\sigma_2\circ\sigma_1$ 
has the positive topological entropy $\log\left(9+4\sqrt{5}\right)$. 
In this case also, there are no $f$-periodic curves. 
we refer the reader to \cite{Bar, Ma, Wa}, 
for its arithmetical properties. See also \cite[Introduction]{Mc}. 
\end{Example}

\begin{Example}[More K3 surface automorphisms]
\label{subsec:eg:2}
Let $S$ be an even lattice of signature $(1,1)$,
which represents neither $0$ nor $-2$. Then
$O(S) := \{g: S \to S \mid \text{$g$ is a lattice isometry}\}$
is infinite, and there is $g_S$ such that the eigenvalues of $g$ are
given by $\lambda, \lambda^{-1}$ for some $\lambda >1$
(cf. \cite[Example~in~p.~581]{PS}).

Let $L_{K3}$ denotes the $K3$ lattice, i.e., the even unimodular lattice 
of signature $(3, 19)$. 
We take a primitive embedding $S \hookrightarrow L_{K3}$ 
(\cite[Chap.~I, Theorem~(2.9)]{BHPV}), 
by which we regard $S$ as a sublattice of $L_{K3}$.
Let $T$ be an orthogonal complement of $S$ in $L_{K3}$.
We set $g := (g_S, 1_T)$ acting on $S \oplus T$.
Since $S^*/S$ is a finite group,
by replacing $g$ by $g^n$ for some positive integer $n$,
we may assume that $g$ acts on $S^*/S$ trivially.
Then $g$ extends to a lattice isometry $g: L_{K3} \to L_{K3}$.

We take a $K3$ surface $X$ and an isometry $\rho: H^2(X, \ZZ) \to
L_{K3}$ such that the N{\'e}ron-Severi lattice and the transcendental
lattice of $X$ are respectively given by $S$ and $T$ via $\rho$. Since
$g$ acts on $T$ identically, $g$ is a Hodge isometry. Moreover, since
$S$ does not represent $-2$, $g$ preserves the ample cone.  Then 
Torelli's theorem (\cite[Chap.~VIII, Proposition~(6.2)]{BHPV}) 
yields that there is an automorphism $f: X \to X$ such 
that $f^*$ is conjugate to $g$ via $\rho$. 
The automorphism $f$ has topological entropy $\log\lambda >0$.
\end{Example}

\begin{Example}[Abelian surfaces and Kummer surfaces]
\label{subsec:eg:3}
Following \cite{Ca} and \cite[\S4]{Mc}, 
we present some automorphisms of abelian surfaces and 
Kummer surfaces of positive topological entropy.  
(Some of the notation here will be used in the proof of 
Lemma~\ref{lemma:D}.)

Let $X$ be an abelian surface over $\CC$, and 
$f: X \to X$ an automorphism. 
We fix an analytic isomorphism $X \simeq \CC^2/\Lambda$, 
where $\Lambda$ is a lattice in $\CC^2$.
We put $v := f(0)$, 
and let $t_v: X \to X$ the translation 
by $v$. Then there is a group isomorphism 
$g: X \to X$ such that 
$f = t_v\circ g$. 

Let $V \in \CC^2$ be a lift of $v$, and 
$T_V: \CC^2 \to \CC^2$ the translation by $V$. 
We denote by $G: \CC^2 \to \CC^2$ the lift of 
$g$ such that $G(0)=0$. 
Then $F:= T_V\circ G$ is the lift of $f$ such that $F(0)= V$. 
We take the matrix $P = \begin{pmatrix}p & q \\ r & s \end{pmatrix} 
\in \GL_2(\CC)$ representing $G$. 
Let $\alpha, \beta$ be the roots of the characteristic polynomial of 
$P$ such that $|\alpha| \geq |\beta|$. 

By identifying $H^{1,1}(X, \RR)$ as the space of translation-invariant 
$(1,1)$-forms on $X$, we compute the action of $f^*$ 
on $H^{1,1}(X, \RR)$ to find that the spectral radius $\lambda$ of 
$\rest{f^*}{H^{1,1}(X, \RR)}$ is given by 
$\lambda = |\alpha|^2$ (cf. \cite[\S4]{Mc}). 
We also note that $|\alpha\beta| = |ps -qr| = 1$ 
(cf. \cite[\S4]{Mc}). 
In particular, if $|p+s| > 2$, then $f$ has positive topological 
entropy. 

Let $[-1]: X \to X$ be the inverse morphism. 
Let $K$ be the Kummer surface associated with $X$, 
i.e., the $K3$ surface obtained by resolving the $16$ 
ordinary double points of $X/[-1]$. 
In the sequel, we assume $v=0$, i.e., $f=g$. In this case, 
since $f$ commutes with $[-1]$, 
$f$ induces the isomorphism $\overline{f}: K \to K$. 
By \cite[\S4]{Mc}, $\overline{f}$ has the same topological 
entropy as $f$. 

Let $C$ be a nodal curve on $K$ that corresponds to one of the 
$16$ ordinary double points of $X/[-1]$. Then $C$ is 
$\overline{f}$-periodic. Specifically, if $C_0$ be the nodal curve 
which corresponds to the origin of 
$X/[-1]$, then $\overline{f}(C_0) = C_0$ and 
$\rest{\overline{f}}{C_0}: C_0 \to C_0$ 
is conjugate to 
\[
\PP^1 \ni (x_0:x_1) \mapsto (px_0+ qx_1: rx_0 + sx_1) \in \PP^1
\] 
via a suitable isomorphism $C \cong \PP^1$. 

For concrete examples, we consider an abelian surface $X$ which is the 
product of an elliptic curve $E$ and itself: $X = E \times E$. 
Then for each $P= \begin{pmatrix}p & q \\ r & s \end{pmatrix} 
\in \SL_2(\ZZ)$ with $|p+s| > 2$, 
the group isomorphism $f: X \to X$ given by 
$f(x, y) = (px + q y, rx+ sy)$ induces 
the automorphism $\overline{f}: K \to K$ of positive topological entropy. 
\end{Example}

In the rest of this section, 
we consider a version of Proposition~\ref{prop:nu} for 
$\RR$-divisors modulo $\RR$-linearly equivalence class.  

\begin{Lemma}
\label{lemma:D}
Let the notation and assumption be as in Proposition~\ref{prop:nu}. 
\begin{enumerate}
\item[(1)]
Then there are $\RR$-divisors $D_+$ and $D_-$ on $X$ with 
$[D_+]=\nu_+$ and  $[D_-]=\nu_-$ that satisfy  
\[
  f^*(D_+) \sim_{\RR} \lambda D_+
  \quad\text{and}\quad 
  f^*(D_-) \sim_{\RR} \lambda^{-1} D_-.
\] 
Further, if $X$ is not birational to an abelian surface, 
then such $D_+$ and $D_-$ are unique up to $\RR$-linearly equivalence 
classes. 
\item[(2)]
We set $D := D_+ + D_-$. Then $D$ is an $\RR$-divisor on $X$ 
with $[D]=\nu$ that satisfies 
\[
  f^*(D) + f^{-1*}(D) \sim_{\RR} (\lambda + \lambda^{-1}) D. 
\]
\end{enumerate} 
\end{Lemma}

\Proof
Since the assertion (2) follows from (1), it suffices to show (1). 

{\bf Case 1:} 
Suppose $X$ is not birational to an abelian surface. Then 
$X$ is either a rational surface, or a surface birational to 
a $K3$ surface or an Enriques surface. In either case, we have 
$H^1(X, \OO_X) =0$. Then the exponential sequence yields 
\[
0 \longrightarrow \Cl(X) \overset{c_1}{\longrightarrow} H^2(X, \ZZ). 
\]
Hence $\Cl(X)_{\RR} \cong N^1(X)_{\RR}$ via $c_1$.  
Thus the existence and uniqueness of $D, D_+, D_-$ follows. 

{\bf Case 2:} 
Suppose $X$ is birational to an abelian surface.  
We will show the existence of $D_+$. The existence of $D_-$ is shown
similarly. 
From the arguments in the proof of Lemma~\ref{lemma:contraction},  
we may assume that $X$ is an abelian surface. 

Since $\lambda$ is an eigenvalue of $f^*: N^1(X) \to N^1(X)$, 
$\lambda$ is a unit, i.e., both $\lambda$ and $\lambda^{-1}$ are 
algebraic integers. We put $k = [\QQ(\lambda):\QQ]$. Let 
\[
  S(t) := t^k + c_1 t^{k-1} + \cdots + c_{k-1} t + c_k
\]
be the irreducible monic polynomial with $S(\lambda)=0$. 
By Theorem~\ref{thm:lambda}(2), the roots of $S(t)$ are given by 
\[
  \{\lambda, \lambda^{-1}, \gamma_1, \overline{\gamma_1}, 
  \ldots, \gamma_{\frac{k}{2}-1}, \overline{\gamma_{\frac{k}{2}-1}}\},
\]
where $|\gamma_i|=1$. 

Since $\lambda$ is an eigenvalue of $f^*$, there is a nef class 
$\nu_+^{\circ} \in N^1(X)\otimes_{\ZZ} \ZZ[\lambda]$ with 
$f^*(\nu_+^{\circ}) = \lambda \nu_+^{\circ}$. 
Since $\nu_+ = c \nu_+^{\circ}$ for some $c >0$, 
the existence of an $\RR$-divisor $D_+^{\circ}$ for 
$\nu_+^{\circ}$ (i.e., $[D_+^{\circ}] = \nu_+^{\circ}$ and 
$f^*(D_+^{\circ}) \sim_{\RR} D_+^{\circ}$) implies the existence 
of an $\RR$-divisor $D_+$ for $\nu_+$. Indeed, we just put 
$D_+ := c D_+^{\circ}$. Hence we may assume that 
$\nu_+ = \nu_+^{\circ}$. 

Since $\Cl(X) \to N^1(X)$ is surjective, we take 
$D_+^{'} \in \Cl(X)\otimes_{\ZZ} \ZZ[\lambda]$ with $[D_+^{'}] = \nu_+$. 
Let $\Cl^0(X)$ the kernel of $\Cl(X) \to N^1(X)$. 
Then $T := f^*(D_+^{'}) - \lambda D_+^{'}$ belongs to 
$\Cl^0(X)\otimes_{\ZZ} \ZZ[\lambda]$. 

\begin{Claim}
\label{claim:lemma:D}
The map $f^* - \lambda\id : \Cl^0(X)\otimes_{\ZZ} \ZZ[\lambda] 
\to \Cl^0(X)\otimes_{\ZZ} \ZZ[\lambda]$ 
is surjective. 
\end{Claim}

We will show the claim at the end. For the moment, we admit the claim 
and proceed the proof. 
Admitting Claim~\ref{claim:lemma:D}, there is 
$T' \in \Cl^0(X)\otimes_{\ZZ}\ZZ[\lambda]$ such 
that $f^*(T') - \lambda T' = T$. 
We put $D_+ := D_+^{'} - T'$. Then $[D_+] = \nu_+$ and 
$f^*(D_+) - \lambda D_+ = 0 \in \Cl(X)\otimes_{\ZZ} \ZZ[\lambda]$. 
Hence we get the existence of $D_+$ for $\nu_+$. 

We are left to prove Claim~\ref{claim:lemma:D}. To see this, 
let us see the action of $f^* - \lambda\id$ on 
$\Cl^0(X)\otimes_{\ZZ}\ZZ[\lambda]$ in more detail. 
We use the notation in Example~\ref{subsec:eg:3}. 

Let $\overline{\Omega}:= \Hom_{\overline{\CC}}(\CC^2, \CC)$ be 
the space of $\CC$-antilinear forms $l: \CC^2 \to \CC$. 
Let $\widehat{\Lambda} := 
\{l \in \overline{\Omega} \mid 
\Image l(\Lambda) \subseteq \ZZ\}$ be the dual lattice of $\Lambda$. 
The dual abelian surface of $X$ is defined by 
$\widehat{X} := \overline{\Omega}/\widehat{\Lambda}$. 
By \cite[Proposition~2.4.1]{BL}, there is a canonical isomorphism 
$\widehat{X} \simeq \Cl^0(X)$. 
We define $\widehat{G}: \overline{\Omega} \to 
\overline{\Omega}$ by $l \mapsto l \circ G$. 
Then $\widehat{G}$ induces an automorphism 
$\widehat{g}: \widehat{X} \to \widehat{X}$. 
Then, by \cite[Eqn.\!~(1)~in~p.~35]{BL}, one sees that 
$f^* = \widehat{g}$ under the canonical isomorphism 
$\Cl^0(X) \simeq \widehat{X}$. (We note that, for $L \in \Cl^0(X)$ and 
$v \in X(\CC)$, we have $t_{v}^*(L) = L$ 
by \cite[Lemma~2.3.2,~Proposition~2.4.1]{BL}.)

We define the basis $\{l_1, l_2\}$ of 
$\overline{\Omega} = \Hom_{\overline{\CC}}(\CC^2, \CC)$ by 
$l_1(z_1, z_2) = \overline{z_1}$ and 
$l_2(z_1, z_2) = \overline{z_2}$. Then $\widehat{G}$ is represented 
by the matrix ${}^t \overline{P} 
= \begin{pmatrix}\overline{p} & \overline{r} \\ \overline{q} 
& \overline{s} \end{pmatrix}$. 

To show the claim, it suffices to show 
\[
  \widehat{G}\otimes_{\ZZ}\ZZ[\lambda] - \lambda\id :
  \overline{\Omega}\otimes_{\ZZ}\ZZ[\lambda] \to
  \overline{\Omega}\otimes_{\ZZ}\ZZ[\lambda] 
\]
is surjective. (Note that $\lambda$ in $\lambda\id$ acts 
on $\overline{\Omega}\otimes_{\ZZ}\ZZ[\lambda]$ through 
$\ZZ[\lambda]$.) We take the $\ZZ$-basis 
$\{1, \lambda, \cdots, \lambda^{k-1}\}$ 
of $\ZZ[\lambda]$ over $\ZZ$. Then, with the basis 
$\{l_1 \otimes 1, l_2\otimes 1, 
l_1 \otimes \lambda, l_2\otimes \lambda, \ldots, 
l_1 \otimes \lambda^{k-1}, l_2\otimes \lambda^{k-1}
\}$ of the $\CC$-vector space 
$\overline{\Omega}\otimes_{\ZZ}\ZZ[\lambda]$, 
$\widehat{G}\otimes_{\ZZ}\ZZ[\lambda] - \lambda\id$ 
is represented by the $(2k, 2k)$ matrix 
\[
  M = 
  \begin{pmatrix}
  {}^t \overline{P} &   0   &  \cdots  &0    & c_k I \\
  -I          &   {}^t \overline{P}    &  \cdots &0     & c_{k-1} I \\
   \vdots     &    \ddots    & \ddots & \vdots  & \vdots\\
  0 & 0 & \ddots & {}^t \overline{P} & c_{2} I\\
  0 & 0 & \cdots &  -I     & {}^t \overline{P} + c_1 I
  \end{pmatrix}
  \in M(2k,2k;\CC). 
\]
We have $\det M = S(\overline{\alpha}) 
S(\overline{\beta})$, 
where $\overline{\alpha}, \overline{\beta}$ are the eigenvalues 
of ${}^t \overline{P}$ with 
$|\overline{\alpha}| = \sqrt{\lambda}$ and 
$|\overline{\beta}| =\sqrt{\lambda^{-1}}$. 

Since $\lambda, \lambda^{-1}$ are the only roots of $S(t)$ 
outside the unit circle, we get 
$S(\overline{\alpha})\neq 0$ and $S(\overline{\beta}) \neq 0$. 
Hence $\det M \neq 0$. It follows that 
$\widehat{G}\otimes_{\ZZ}\ZZ[\lambda] - \lambda\id$ is 
surjective. 
\QED

\begin{Remark}
\label{remark:quadratic:2}
If $\lambda$ is a quadratic irrational number, 
then, with a suitable $\nu$, $D$ can be taken as 
an integral divisor. (cf. Remark~\ref{remark:quadratic}). 
\end{Remark}

\medskip
\setcounter{equation}{0}
\section{Preliminaries on height functions}
\label{sec:height}
In the latter half of this paper, we turn our attention to 
arithmetic properties. In this section, we briefly review some 
facts on height functions and fix notation. 
A general reference is \cite[Part~B]{HS} for example. 

Let $K$ be a number field, and $O_K$ its ring of integers. 
For $x = (x_0:\cdots:x_N) \in \PP^N(K)$, the logarithmic naive 
height of $x$ is defined by 
\[
  h_{nv}(x) = 
  \frac{1}{[K:\QQ]} 
  \left[
  \sum_{P \in \Spec(O_K)\setminus\{0\}} 
  \log\max_{i}\{\Vert x_i\Vert_{P}\} 
  + \sum_{\sigma: K \hookrightarrow \CC} 
  \log\max_{i}\{\vert \sigma(x_i)\vert\} 
  \right], 
\]
where $\Vert x_i\Vert_P = \#(O_K/P)^{-\ord_{P}(x_i)}$. 
This definition naturally extends to all points 
$x \in \PP^N(\overline{\QQ})$, 
and gives rise to the logarithmic naive height function 
$h_{nv}: \PP^N(\overline{\QQ}) \to \RR$. 

We state Weil's height machine, which associates a height function 
to every $\RR$-divisor, and summarize some properties of height functions. 
By a divisor we mean a Cartier divisor. 
For a proof, we refer the reader to \cite[Theorem~B.3.2]{HS}. 
See also \cite[Theorems~1.1.1,~1.1.2]{Ka}.  
In the following, $O(1)$ denotes a bounded function on 
$X(\overline{K})$. 

\begin{Theorem}[cf. \cite{HS} Theorem~B.3.2]
\label{thm:height:machine}
Let $K$ be a number field. 
For any projective variety $X$ over $K$ 
and any $\RR$-divisor $D$ on $X$, 
there is a way to attach a real-valued function,  
called a height function corresponding to $D$, 
\[
  h_{X, D}: X(\overline{K}) \to \RR
\]
with the following properties. 
\begin{enumerate}
\item[(i)]
\textup{(Additivity)} 
$h_{X, D_1 + D_2} = h_{X, D_1} + h_{X, D_2} + O(1)$ for any 
$D_1, D_2 \in \Div_{\RR}(X)$. 
\item[(ii)]
\textup{(Normalization)} 
If $H$ is a hyperplane in $\PP^N$,  
then $h_{\PP^N, H} = h_{nv} + O(1)$. 
\item[(iii)]
\textup{(Functoriality)} 
If $f: X \to Y$ is a morphism of projective varieties 
and $D \in  \Div_{\RR}(Y)$ with $f(X) \not\subseteq \Supp(D)$, 
then $h_{X, f^*D} = h_{Y, D} \circ f + O(1)$.  
\end{enumerate}
The height functions $h_{X,D}$ are uniquely determined,  
up to $O(1)$, by the above three properties.  
Moreover, they satisfy the following 
properties. 
\begin{enumerate}
\item[(iv)]
\textup{(Linear equivalence)}  
If $D_1, D_2 \in \Div_{\RR}(X)$ are $\RR$-linearly equivalent, 
then $h_{X, D_1} = h_{X, D_2} + O(1)$. 
\item[(v)]
\textup{(Northcott finiteness property)} 
Assume $D \in \Div_{\RR}(X)$ is ample. 
Then for any positive integer $d$ and 
real number $M$, the set 
\[
  \{x \in X(\overline{K}) \mid 
  [K(x): K] \leq d, \; h_{X, D}(x) \leq M\}
\]
is finite.  
\item[(vi)]
\textup{(Positivity)}
If $D \in \Div_{\RR}(X)$ is effective, 
then there is a constant $c$ such that 
$h_{X, D}(x) \geq c$ 
for all $x \in (X\setminus \Supp(D))(\overline{K})$. 
If $D' \in \Div_{\RR}(X)$ is ample, 
then there is a constant $c'$ such that 
$h_{X, D'}(x) \geq c'$
for all $x \in X(\overline{K})$. 
\item[(vii)]
If $D_1, D_2 \in \Div_{\RR}(X)$ are ample, 
then there are positive constants $a_1, a_2 >0$ and constants 
$b_1, b_2 \in\RR$ such that 
\[
a_1 h_{X, D_1}(x) + b_1 \leq
h_{X, D_2}(x) 
\leq a_2 h_{X, D_1}(x) + b_2
\]
for all $x \in X(\overline{K})$. 
\end{enumerate}
\end{Theorem}

A function that coincides with $h_{X, D}$ up to $O(1)$ is also 
called a height function corresponding to $D$. In virtue of 
Theorem~\ref{thm:height:machine}(iv), 
for a linear equivalence class $D \in \Cl(X)\otimes\RR$, 
we define $h_{X, D}$ to be a height function corresponding 
an $\RR$-divisor which represents $D$. 
We will often write $h_D$ in place of $h_{X, D}$ 
when $X$ is obvious from the context. 

\medskip
\setcounter{equation}{0}
\section{Canonical height functions}
\label{sec:canonical}
Let $X$ be a smooth projective surface over a number field 
$K \subset \CC$, and $f: X \to X$ an automorphism over $K$.
Let $f_{\CC}: X_{\CC} \to X_{\CC}$ be the automorphism induced 
from $f$, where  $f_{\CC} := f \times_{\Spec(K)} \Spec(\CC)$ and 
$X_{\CC} := X \times_{\Spec(K)} \Spec(\CC)$.  
We define the topological entropy of $f$ by that of $f_{\CC}$. 

\begin{Lemma}
\label{lemma:field:def}
The $\RR$-divisor classes $\nu, \nu_+, \nu_- \in N^1(X_{\CC})_{\RR}$ in 
Proposition~\ref{prop:nu}, the maximal $f$-invariant curve $E$ on $X_{\CC}$ 
in Eqn.\!~\eqref{eqn:maximal:curve}, and the $\RR$-divisors $D, D_+, D_-$ 
on $X_{\CC}$ in Lemma~\ref{lemma:D} are all defined 
over $\overline{K}$. 
\end{Lemma}

\Proof
We set $X_{\overline{K}} := X \times_{\Spec(K)} \Spec(\overline{K})$. 
We claim that $N^1(X_{\overline{K}}) = N^1(X_{\CC})$ naturally. 
Noting that a divisor refers to a Cartier divisor, $\Cl(S)$ is 
canonically isomorphic to $\Pic(S)$ for a projective scheme $S$  
over a field. In this proof, we use $\Pic$ instead of $\Cl$. 

Let $\Pic_{X_{\overline{K}}/\overline{K}}$ be the Picard 
scheme, and $\Pic^0_{X_{\overline{K}}/\overline{K}}$ the connected 
component of $\Pic_{X_{\overline{K}}/\overline{K}}$ that contains 
the trivial line bundle. In our case, 
$\Pic^0_{X_{\overline{K}}/\overline{K}}$ is also 
the irreducible component of $\Pic_{X_{\overline{K}}/\overline{K}}$ 
that contains the trivial line bundle.  

By \cite[\S8.1~Proposition~4]{BLR}, 
we have $\Pic_{X_{\overline{K}}/\overline{K}}(\overline{K}) 
= \Pic(X_{\overline{K}})$ and 
$\Pic_{X_{\overline{K}}/\overline{K}}(\CC) 
= \Pic(X_{\CC})$. 
Since, for an integral domain $A$ 
over $\overline{K}$, $A \otimes_{\overline{K}} \CC$ is an integral 
domain, we have  $\Pic^0_{X_{\overline{K}}/\overline{K}}(\overline{K}) 
= \Pic^0(X_{\overline{K}})$ and $\Pic^0_{X_{\overline{K}}/\overline{K}}(\CC) 
= \Pic^0(X_{\CC})$. 
Hence $\Pic(X_{\overline{K}})/\Pic^0(X_{\overline{K}}) 
= \Pic(X_{\CC})/\Pic^0(X_{\CC}) \ 
(= \Pic_{X_{\overline{K}}/\overline{K}}/
\Pic^0_{X_{\overline{K}}/\overline{K}}$). 

By Matsusaka's theorem \cite{Matsu}, we have 
$N^1(X_{\overline{K}}) 
= \left(\Pic(X_{\overline{K}})/\Pic^0(X_{\overline{K}})\right)/ 
(\text{torsion})$  and $ N^1(X_{\CC}) = 
\left(\Pic(X_{\CC})/\Pic^0(X_{\CC})\right) / (\text{torsion})$.  
Hence $N^1(X_{\overline{K}}) = N^1(X_{\CC})$. 
We then see that $\nu, \nu_+, \nu_-$ in 
Proposition~\ref{prop:nu} are defined over $\overline{K}$. 

Next, we show that $E$ is defined over $\overline{K}$. 
We fix an ample line bundle $\OO(1)$ on $X_{\overline{K}}$ and let 
$P$ be the Hilbert polynomial of $E$. Let $\calH := 
\Hilb_P(X_{\overline{K}}/\overline{K})$ be the Hilbert scheme. 
Let $\overline{e}: \Spec(\CC) \to \calH$ 
be the $\CC$-valued point corresponding to $E$, 
and $e$ the scheme point, i.e., the image of $\overline{e}$ 
in $\calH$. Let $\kappa$ be the residue field of $\OO_{\calH, e}$. 
Via $\Spec(\CC) \overset{\overline{e}}{\rightarrow} \calH 
\to \Spec(\overline{K})$, we have $\overline{K} \subset \kappa \subset \CC$. 
We claim that $\kappa = \overline{K}$. To see this, suppose 
$\kappa \neq \overline{K}$. Then there is an automorphism 
$\sigma: \CC \to \CC$ such that $\rest{\sigma}{\overline{K}} = \id$ 
and such that $\sigma(\kappa) \neq \kappa$. 
Let $\sigma^*: \Spec(\CC) \to \Spec(\CC)$ 
and $\sigma^*: X_{\CC} \to X_{\CC}$ be morphisms induced 
by $\sigma$. Since $E$ is the (unique) maximal $f$-invariant curve, 
we see that $\sigma^*(E) = E$. 
Hence $\overline{e} \circ \sigma^* =  \overline{e}$, 
and $\sigma(\kappa) = \kappa$. This is a contradiction. 
Thus $\kappa = \overline{K}$, and $\overline{e}$ decomposes into 
$\Spec(\CC) \to \Spec(\overline{K}) \to \calH$. 
Hence $E$ is defined over $\overline{K}$. 

Suppose $X$ is not birational to an abelian surface. Then 
$\Cl(X_{\overline{K}}) \simeq N^1(X_{\overline{K}})$ implies that 
$D, D_+, D_-$ in Lemma~\ref{lemma:D} are defined over $\overline{K}$.

Suppose now $X$ is birational to an abelian surface. 
We see that 
the morphism $f^* - \lambda\id : 
\Pic^0_{X_{\overline{K}}/\overline{K}} \otimes_{\ZZ} \ZZ[\lambda] \to 
\Pic^0_{X_{\overline{K}}/\overline{K}} \otimes_{\ZZ} \ZZ[\lambda] $ 
is surjective, because 
$f^* - \lambda\id: \Pic^0_{X_{\overline{K}}/\overline{K}}(\CC) 
\otimes_{\ZZ} \ZZ[\lambda] \to
\Pic^0_{X_{\overline{K}}/\overline{K}}(\CC)
\otimes_{\ZZ} \ZZ[\lambda]$ is surjective 
by Claim~\ref{claim:lemma:D}. Hence 
$f^* - \lambda\id: 
\Pic^0(X_{\overline{K}}) \otimes_{\ZZ} \ZZ[\lambda] 
\to \Pic^0(X_{\overline{K}}) \otimes_{\ZZ} \ZZ[\lambda]$ 
is surjective. A similar argument as in the proof of 
Lemma~\ref{lemma:D} shows that $D, D_+, D_-$ 
are defined over $\overline{K}$.
\QED

In virtue of Lemma~\ref{lemma:field:def}, we will often replace $K$ 
by its suitable finite extension so that $D$ and $E$ are defined 
over $K$. 

\begin{Theorem}
\label{thm:canonical:height}
Let $X$ be a smooth projective surface over a number field $K \subset \CC$, 
and $f: X \to X$ an automorphism of positive 
topological entropy $\log\lambda >0$. 
Let $D$ be a nef and big $\RR$-divisor in Lemma~\ref{lemma:D}, 
and $E$ the maximal $f$-invariant curve. Replacing $K$ by a suitable 
finite extension if necessary, we assume that $D$ and $E$ are defined over 
$K$. Then there exists a unique function 
\[
  \widehat{h}_{D}: X(\overline{K}) \to \RR
\]
with the following two properties. 
\begin{enumerate}
\item[(i)]
$\widehat{h}_D$ is a height function corresponding 
to $D$. 
\item[(ii)]
$\widehat{h}_D(f(x)) + \widehat{h}_D(f^{-1}(x))
= \left(\lambda + \lambda^{-1}\right) \widehat{h}_D(x)$ 
for all $x \in X(\overline{K})$.  
\end{enumerate} 
Moreover, $\widehat{h}_{D}$ satisfies the following properties. 
\begin{enumerate}
\item[(iii)] 
$\widehat{h}_D(x) \geq 0$ for all $x \in X(\overline{K})$. 
\item[(iv)]
If $x \in E(\overline{K})$, then $\widehat{h}_D(x) = 0$. 
\item[(v)]
Suppose $x \in (X\setminus E)(\overline{K})$. Then 
$\widehat{h}_D(x) = 0$ if and only if $x$ is $f$-periodic. 
\item[(vi)]
$\widehat{h}_D$ satisfies the Northcott finiteness property 
on $(X\setminus E)(\overline{K})$: For any positive integer $d$ 
and real number $M$, the set
\[
  \{x \in (X\setminus E)(\overline{K})
  \mid [K(x):K] \leq d, \; \widehat{h}_D(x) \leq M \}
\]
is finite. 
\end{enumerate} 
We call $\widehat{h}_{D}$ a 
{\em canonical height function for $f$ (and $D$)}. 
\end{Theorem}

\Proof
(i) and (ii)  This follows from 
a similar argument as in \cite[Theorem~1.1]{SiK3}. 
Since $\lambda >1$ and $f^*(D) + f^{-1 *}(D) \sim_{\RR} 
\left(\lambda + \lambda^{-1}\right) D$ by 
Lemma~\ref{lemma:D}(2), 
Theorem~\ref{thm:canonical:height}(i)(ii) is also 
a special case of \cite[Theorem~1.2.1]{Ka}, once we regard $(X; f, f^{-1})$ 
as a dynamical eigensystem for $D$ of degree $(\lambda + \lambda^{-1})$. 

For the sake of completeness, we briefly sketch a proof 
in line with \cite{SiK3}. 
Let $D_+, D_-$ be as in Lemma~\ref{lemma:D}, 
and $h_{D_+}, h_{D-}$ any height functions corresponding 
to $D_+, D_-$, respectively. Since 
\begin{equation}
\label{eqn:canonical:height:0}
  f^*(D_+) \sim_{\RR} \lambda D_+
  \quad\text{and}\quad 
  f^*(D_-) \sim_{\RR} \lambda^{-1} D_-, 
\end{equation}
Tate's telescoping argument asserts that the following limits exist: 
\begin{equation}
\label{eqn:canonical:height:a}
  \widehat{h}_{D_+}(x) := 
  \lim_{n\to +\infty} \frac{1}{\lambda^n} h_{D_+}\left(f^n(x)\right)  
  \quad\text{and}\quad 
  \widehat{h}_{D_-}(x) :=  
  \lim_{n\to +\infty} \frac{1}{\lambda^n} h_{D-}\left(f^{-n}(x)\right).   
\end{equation}
Then $\widehat{h}_{D_+}$ and $\widehat{h}_{D_-}$ are height functions 
corresponding to $D_+$ and $D_-$, which satisfy 
$\widehat{h}_{D_+} \circ f = \lambda \widehat{h}_{D_+}$ and   
$\widehat{h}_{D_-}\circ f = \lambda^{-1} \widehat{h}_{D_-}$, 
respectively. One defines
\begin{equation}
  \label{eqn:canonical:height:b}
  \widehat{h}_D := \widehat{h}_{D_+} + \widehat{h}_{D_-}. 
\end{equation}
One can check that $\widehat{h}_D$ satisfies (i) and (ii). 
If $\widehat{h}_D$ and $\widehat{h}^{\prime}_D$ are 
two such functions, then $\widehat{h}_D - \widehat{h}^{\prime}_D$ 
is bounded and satisfies
$(\lambda + \lambda^{-1})\, \Vert \widehat{h}_D 
- \widehat{h}^{\prime}_D\Vert_{\sup} 
\leq 2\, \Vert \widehat{h}_D - \widehat{h}^{\prime}_D\Vert_{\sup}$,  
from which the uniqueness follows. 

We remark that it follows from 
\eqref{eqn:canonical:height:0} that 
\begin{equation}
\label{eqn:canonical:height:c}
  f^{n*}(D) + f^{-n *}(D) \sim_{\RR} 
  \left(\lambda^n + \lambda^{-n}\right) D
\end{equation}
for any $n \in \ZZ$. It follows from 
$\widehat{h}_{D_{\pm}} \circ f = \lambda^{\pm} \widehat{h}_{D_{\pm}}$ 
and \eqref{eqn:canonical:height:b} that 
\begin{equation}
\label{eqn:canonical:height:2}
  \widehat{h}_D(f^n(x)) + \widehat{h}_D(f^{-n}(x))
  = \left(\lambda^n + \lambda^{-n}\right) \widehat{h}_D(x)
\end{equation}
for any $n \in \ZZ$ and $x \in X(\overline{K})$. 
In the following, for clarity of the argument, we prove 
Theorem~\ref{thm:canonical:height}(iii)--(vi) 
in the order of (iv), (iii), (vi), (v). 

\smallskip
(iv) 
To show (iv), we may replace $K$ by a suitable finite extension 
so that every component of $E$ is geometrically irreducible. 
Let $x \in E(\overline{K})$. Let $C$ be an irreducible component 
of $E$ such that $x \in C(\overline{K})$. Since $C$ is an $f$-periodic 
curve, there is a positive integer $n$ such that $f^n(C) = C$. 
We take $D' \sim_{\RR} D$ such that any component of $\Supp(D')$ is 
not equal to $\Supp(C)$, and we set $L := \rest{D'}{C}$. 
We take the normalization $\varphi: \widetilde{C} \to C$, 
set $\widetilde{L} := \varphi^*(L)$, and write 
$\widetilde{f}^n: \widetilde{C} \to \widetilde{C}$ for 
the the automorphism induced by $f^n$.  
Then $\widetilde{L}$ is an $\RR$-divisor on $C$ such that
$\widetilde{f}^{n*}(\widetilde{L}) +
\widetilde{f}^{-n*}(\widetilde{L}) \sim_{\RR} (\lambda^n +
\lambda^{-n}) \widetilde{L}$. It follows from 
Lemma~\ref{lemma:canonical:height} below that 
$\widetilde{L}\sim_{\RR} 0$. 

We put $\widetilde{h} := \left(\rest{\widehat{h}_{D}}{C}\right)
\circ\varphi$. Since $\widetilde{h}$ is a height function corresponding 
to $\widetilde{L}\sim_{\RR} 0$, $\widetilde{h}$ is a bounded function. 
Further, since $\widetilde{h}$ satisfies Eqn.\!~\eqref{eqn:canonical:height:2}
with $\widetilde{h}$ and $\widetilde{f}^n$ 
in place of $\widehat{h}_D$ and $f^n$, 
we conclude $\widetilde{h} \equiv 0$. Thus $\widehat{h}_{D}(x) =0$. 

\smallskip
(iii)
By (iv), we have only to show that $\widehat{h}_{D}(x) \geq 0$ 
for all $x \in (X\setminus E)(\overline{K})$. 
By Proposition~\ref{prop:nu:2} and Proposition~\ref{proposition:finiteness}, 
there is an effective divisor $Z$ with 
$\Supp(Z) \subseteq \Supp(E)$ such that 
$D - \varepsilon Z$ is ample for sufficiently small 
$\varepsilon >0$. Let $h_Z$ be a height function corresponding 
to $Z$. Then by Theorem~\ref{thm:height:machine}(vi), 
there is constant $c_1$ such that 
\begin{equation}
\label{eqn:canonical:height:3}
  \widehat{h}_D(x) - \varepsilon h_Z(x) \geq c_1
  \qquad\text{for all $x \in X(\overline{K})$}. 
\end{equation}
Since $\Supp(Z) \subseteq \Supp(E)$,  
there is constant $c_2$ such that 
\begin{equation}
\label{eqn:canonical:height:4}
  h_Z(x) \geq c_2
  \qquad\text{for all $x \in (X\setminus E)(\overline{K})$}. 
\end{equation}
Note that $f$ induces the automorphism of $X\setminus E$. 
Then Eqn.\!~\eqref{eqn:canonical:height:2} yields that, 
for any $x \in (X\setminus E)(\overline{K})$, 
\[
  \widehat{h}_D(x) = 
  \frac{1}{\lambda^n + \lambda^{-n}} \left(
  \widehat{h}_D(f^n(x)) + \widehat{h}_D(f^{-n}(x))\right)
  \geq 
  \frac{1}{\lambda^n + \lambda^{-n}} 2(c_1 +\varepsilon c_2). 
\]
Letting $n \to +\infty$, we get $\widehat{h}_{D}(x) \geq 0$ 
for $x \in (X\setminus E)(\overline{K})$. 

\smallskip
(vi)
If $x \in (X\setminus E)(\overline{K})$ satisfies $\widehat{h}_D(x) \leq M$, 
then $\widehat{h}_D(x) - \varepsilon h_Z(x) \leq M - \varepsilon c_2$ 
by Eqn.\!~\eqref{eqn:canonical:height:4}. Hence 
\begin{multline*}
\{x \in (X\setminus E)(\overline{K})
\mid [K(x):K] \leq d, \; \widehat{h}_D(x) \leq M \} \\
\quad\subseteq\quad
\{x \in X(\overline{K})
\mid [K(x):K] \leq d, \; \widehat{h}_D(x) - \varepsilon h_Z(x) 
\leq M - \varepsilon c_2\}. 
\end{multline*}
Since $\widehat{h}_D - \varepsilon h_Z$ is a height function 
corresponding to an ample $\RR$-divisor, the latter set is a finite set by 
Theorem~\ref{thm:height:machine}(v). 

\smallskip
(v)
Assume $x$ is $f$-periodic. Then $f^n(x) = x$ for 
some positive integer $n$. By Eqn.\!~\eqref{eqn:canonical:height:2}, 
we get $\widehat{h}_D(x) = 0$. 
On the other hand, suppose $x \in (X\setminus E)(\overline{K})$ 
satisfies $\widehat{h}_D(x) = 0$. Then by (ii) and (iii), we get 
$\widehat{h}_D(f^n(x)) = 0$ for all $n \in \ZZ$. 
By (vi), we find that 
$\{f^n(x) \mid n\in \ZZ\}$ is a finite set. 
It follows that $x$ is an $f$-periodic point. 
\QED

\begin{Lemma}
\label{lemma:canonical:height}
Let $C$ be a smooth projective curve, 
$f$ an automorphism of $C$, and $\lambda >1$ 
a positive number. Suppose that 
an $\RR$-divisor $L$ on $C$ satisfies
\begin{equation} 
\label{eqn:canonical:height:c:similar}
  f^*(L) + f^{-1*}(L) \sim_{\RR} (\lambda + \lambda^{-1}) L.
\end{equation}   
Then $L \sim_{\RR} 0$. 
\end{Lemma}

\Proof 
By induction on $n$, we have  
$f^{n*}(L) + f^{-n*}(L) \sim_{\RR} (\lambda^n + \lambda^{-n}) L$, 
so that we may freely 
replace $\lambda$ by $\lambda^n$ and 
$f$ by $f^n$ to prove the lemma. 

Suppose that $g(C) = 0$. Since $f^*(L) \sim_{\RR} L$ and 
$f^{-1*}(L) \sim_{\RR} L$, 
Eqn.\!~\eqref{eqn:canonical:height:c:similar}
implies that $\left(\lambda + \lambda^{-1} -2\right) L \sim_{\RR} 0$. 
Hence $L \sim_{\RR} 0$. 

Next, suppose that $g(C) = 1$. We fix any point $O$ of $C$ and regard $C$ 
as an elliptic curve with the origin $O$. 
Since the set of the group isomorphisms of $C$ is a finite set, 
by replacing $f$ by a suitable $f^n$, we may assume that 
$f$ is a translation. 
By the theorem of the square, we then get $f^*(L) + f^{-1*}(L) 
\sim_{\RR} 2 L$. Thus $L \sim_{\RR} 0$ from 
Eqn.\!~\eqref{eqn:canonical:height:c:similar}. 

Finally, suppose that $g(C) \geq 2$. Replacing $f$ by 
a suitable $f^n$, we may assume that $f = \id$. 
Then Eqn.\!~\eqref{eqn:canonical:height:c:similar}
implies $L \sim_{\RR} 0$. 
\QED

\begin{Remark}
In Example~\ref{subsec:eg:1}, let $\widehat{h}_{\rom{Sil}}$ denote  
the canonical height function on $X$ over a number field, 
which is constructed in \cite[Theorem~1.1]{SiK3}.  
As is checked in \cite[Proposition~1.4.1]{Ka}, 
$\widehat{h}_{\rom{Sil}} = (1+\sqrt{3})\, \widehat{h}_{H_1 + H_2}$. 
\end{Remark}

We need in \S\ref{subsec:non:periodic} 
the following proposition. 

\begin{Proposition}
\label{prop:decomp}
Let the notation be as in Theorem~\ref{thm:canonical:height}. 
Let $\widehat{h}_{D_+}$ and $\widehat{h}_{D_-}$ as in 
Eqn.\!~\eqref{eqn:canonical:height:a}.  
\begin{enumerate}
\item[(1)]
$\widehat{h}_{D_+}(x) \geq  0$ and $\widehat{h}_{D_-}(x) \geq 0$ 
for all $x \in X(\overline{K})$. 
\item[(2)]
For $x \in (X\setminus E)(\overline{K})$, 
one has $\widehat{h}_{D_+}(x) = 0 
\Longleftrightarrow
\widehat{h}_{D_-}(x) = 0 
\Longleftrightarrow
\widehat{h}_{D}(x) = 0$.  
\end{enumerate}
\end{Proposition}

\Proof
We follow \cite{SiK3}.

(1) Since $\widehat{h}_D$ is non-negative by 
Theorem~\ref{thm:canonical:height}(iii), we have, 
for any $x \in X(\overline{K})$, 
\[
\widehat{h}_{D_+}(x) 
= \lambda^{-n} \widehat{h}_{D_+}(f^n(x)) 
\geq - \lambda^{-n} \widehat{h}_{D_-} (f^n(x)) 
= - \lambda^{-2n} \widehat{h}_{D_-} (x).
\]
Letting $n \to +\infty$, we get 
$\widehat{h}_{D_+}(x) \geq 0$. 
Similarly $\widehat{h}_{D_-}(x) \geq 0$. 

(2) Let $L$ be an finite extension field of $K$ over which 
$x$ is defined. Suppose $\widehat{h}_{D_+}(x) = 0$. 
Then we have 
\[
  \widehat{h}_D(f^n(x)) 
  = \widehat{h}_{D_+}(f^n(x)) + \widehat{h}_{D_-}(f^n(x))
  = \lambda^{-n} \widehat{h}_{D_-}(x). 
\]
Thus
\[
  \{f^n(x) \mid n =1, 2, \ldots, \}
  \subseteq 
  \{
  z \in (X\setminus E)(L)
  \mid \widehat{h}_D(z) \leq \widehat{h}_{D_-}(x)\}. 
\]
By Theorem~\ref{thm:canonical:height}(vi), 
$\{f^n(x) \mid n =1, 2, \ldots, \}$ is a finite set. 
It follows that $x$ is $f$-periodic. 
Then Theorem~\ref{thm:canonical:height}(v) yields 
$\widehat{h}_{D}(x) = 0$, so that 
$\widehat{h}_{D_-}(x) = 0$. 
The other directions are shown similarly. 
\QED

\medskip
\setcounter{equation}{0}
\section{Arithmetic properties of $f$}
\label{sec:thm:B:C}
\subsection{Automorphisms of curves}
\label{subsec:p1}
In this subsection, we show some elementary 
arithmetic properties of automorphisms of curves, 
which will be needed in the sequel. 

Let $f: \PP^1 \to \PP^1$ be an automorphism over a number 
field $K$. Then there is a matrix $F = 
\begin{pmatrix} p & q \\ r & s \end{pmatrix} \in 
\GL_2(K)$ such that $f$ is given by 
$f(x:y) = (px+qy:rx+sy)$ for $(x:y)\in \PP^1(\overline{K})$. 
Let $P \in \GL_2(\overline{K})$ 
be a matrix such that $G : = P^{-1} F P \in \GL_2(\overline{K})$ 
is the Jordan canonical form of $F$.  
Then $G$ is one of the following types: 
\begin{enumerate}
\item[(i)]
$G = \begin{pmatrix} \alpha & 0 \\ 0 & \beta
\end{pmatrix}$, where $\frac{\alpha}{\beta}$ is a root of unity; 
\item[(ii)]
$G = \begin{pmatrix} \alpha & 0 \\ 0 & \beta
\end{pmatrix}$, where $\frac{\alpha}{\beta}$ is not a root of unity;  
\item[(iii)]
$G = \begin{pmatrix} \alpha & 1 \\ 0 & \alpha 
\end{pmatrix}$ with $\alpha \in \overline{K}\setminus\{0\}$. 
\end{enumerate}
Let $p: \PP^1 \to \PP^1$ be the automorphism induced from $P$, 
and $g: \PP^1 \to \PP^1$ the automorphism induced from $G$. 

As in \S\ref{sec:height}, 
let $h_{nv}: \PP^1(\overline{K}) \to \RR$ denote the 
logarithmic naive height function. 

First we consider $f$-periodic points.

\begin{Lemma}
\label{lemma:p1}
\begin{enumerate}
\item[(1)]
If $G$ is of type \textup{(i)}, 
then $f$ is periodic, i.e., $f^n = \id$ for some $n \geq 1$. 
\item[(2)]
If $G$ is of type \textup{(ii)}, 
then the set of $f$-periodic points in $\PP^1(\overline{K})$
consists of two points. 
\item[(3)]
If $G$ is of type \textup{(iii)}, 
then the set $f$-periodic points in $\PP^1(\overline{K})$
consists of one point. 
\end{enumerate}
\end{Lemma}

\Proof
Since $f$ is conjugate to $g$, it suffices to show the lemma for 
$g$ in place of $f$. For $g$, the assertions are readily checked. 
\QED

Next we consider non $f$-periodic points. 
We take a non $f$-periodic point $x \in \PP^1(\overline{K})$. 
We write $O_f(x) = \{f^n(x) \mid n \in \ZZ\} 
\subset \PP^1(\overline{K})$ for the orbit 
of $x$ under $f$. 
For $T > 0$, we set 
\[
  N_{f, x}(T) := \# \{y \in O_f(x) \mid h_{nv}(y) \leq T\}, 
\]
which depends on $T$, $O_f(x)$ and $h_{nv}$. 

For real-valued functions $k_1(T), k_2(T)$, we write 
\[
  k_1(T) \gg\ll k_2(T) \quad\textup{as}\quad T \to +\infty 
\]
if there are positive constants 
$a_1, a_2 >0$ such that $a_1 k_1(T) \leq k_2(T) \leq a_2 k_1(T)$ for all 
sufficiently large $T$. As usual, we set $\log^+ x := \log \max\{x, 1\}$ 
for a positive real number $x$. 

\begin{Lemma}
\label{lemma:p2}
\begin{enumerate}
\item[(1)]
If $G$ is of type \textup{(ii)}, then   
$N_{f, x}(T) \gg\ll T$ as $T\to +\infty$. 
\item[(2)]
If $G$ is of type \textup{(iii)}, then  
$N_{f, x}(T) \gg\ll \exp(T)$ as $T\to +\infty$.
\end{enumerate}
\end{Lemma}

\Proof
To prove the lemma, we may replace $K$ by a finite extension 
field. So we may assume that $g: \PP^1 \to \PP^1$ 
and $p: \PP^1 \to \PP^1$ are defined over $K$. 
We take a finite extension field $L$ of $K$ 
over which $x$ is defined. We put $y = (y_0:y_1) 
:= p^{-1}(x) \in \PP^1(L)$.  

By Theorem~\ref{thm:height:machine}(iv), there is a constant $c$ 
such that $\Vert h_{nv} - h_{nv}\circ p\Vert_{\sup} \leq c$ 
on $\PP^1(\overline{K})$. Then 
\[
  \left\vert h_{nv}(f^n(x)) - h_{nv}(g^n(y)) \right\vert 
  =  
  \left\vert h_{nv}(p(g^n (p^{-1}(x))) - h_{nv}(g^n(y)) 
  \right\vert \leq c 
\]
for all $n \in \ZZ$. Hence we get 
$N_{g, y}(T -c) \leq N_{f,x}(T) \leq N_{g, y}(T +c)$. 
This estimate shows that it suffices to show the lemma 
for $g$ and $y$.  

(1) 
We note $y_0y_1\neq 0$. 
Since $\frac{\alpha}{\beta}$ is not a root of unity, 
$h_{nv}((\alpha: \beta)) \neq 0$. 
Since $(\alpha^n : \beta^n) = \left(
(\frac{\alpha}{\beta})^n : 1\right)$
and $g^n(y) = (\alpha^n y_0 : \beta^n y_1) = 
\left((\frac{\alpha}{\beta})^n \frac{y_0}{y_1} : 1\right)$, 
we have  
{\allowdisplaybreaks
  \begin{align}
  \label{eqn:p2:1}
  h_{nv}((\alpha^n : \beta^n)) & = 
  \frac{1}{[L:\QQ]} 
  \left[
  \sum_{P \in \Spec(O_L)\setminus\{0\}} 
  \log^+\left\Vert\left(\frac{\alpha}{\beta}\right)^n\right\Vert_{P} 
  + \sum_{\sigma: L \hookrightarrow \CC} 
  \log^+ \left\vert \sigma\left(\left(
  \frac{\alpha}{\beta}\right)^n\right)\right\vert \right], 
  \\ 
  \label{eqn:p2:2} 
  h_{nv}(g^n(y)) & = 
  \frac{1}{[L:\QQ]} 
  \left[
  \sum_{P \in \Spec(O_L)\setminus\{0\}} 
  \log^+ \left\Vert \left(\frac{\alpha}{\beta}\right)^n \frac{y_0}{y_1} 
  \right\Vert_{P} 
  + \sum_{\sigma: L \hookrightarrow \CC} 
  \log^+ \left\vert \sigma\left(\left(\frac{\alpha}{\beta}\right)^n 
  \frac{y_0}{y_1}\right)\right\vert
  \right]. 
  \end{align}
}

For any real numbers $a, b > 0$, we have 
$\max\{ab, 1\} \leq \max\{a, 1\}\max\{b,1\}$.  
Putting $a = \left(\frac{\alpha}{\beta}\right)^n \frac{y_0}{y_1}$ and 
$b = \left(\frac{y_0}{y_1}\right)^{-1}$ in Eqn.\!~\eqref{eqn:p2:1}, 
and putting $a = \left(\frac{\alpha}{\beta}\right)^n$ 
and $b = \frac{y_0}{y_1}$ in Eqn.\!~\eqref{eqn:p2:2},
we get 
\[
  h_{nv}(g^n(y)) - h_{nv}(y) 
  \leq h_{nv}((\alpha^n: \beta^n)) = 
  |n| h_{nv}((\alpha: \beta))
  \leq h_{nv}(g^n(y)) + h_{nv}(y). 
\]
Then $\lim_{T \to +\infty} 
  \frac{N_{g, y}(T)}{T} = \frac{2}{h_{nv}((\alpha: \beta))}$. 
Hence $N_{g, y}(T) \gg\ll T$ as $T \to +\infty$. 

(2) 
We note $y_1 \neq 0$. 
Since $\frac{1}{\alpha}F \in \GL_2(K)$ is also a lift of 
$f:\PP^1 \to \PP^1$, we may assume that $\alpha=1$.
Then $g^n(y) = (y_0 + n y_1: y_1) = (\frac{y_0}{y_1} + n : 1)$. 
We get 
\[
  h_{nv}(g^n(y)) = 
  \frac{1}{[L:\QQ]} 
  \left[
  \sum_{P \in \Spec(O_L)\setminus\{0\}} 
  \log^+\left\Vert \frac{y_0}{y_1} + n \right\Vert_{P} 
  + \sum_{\sigma: L \hookrightarrow \CC} 
  \log^+ \left\vert \sigma\left(
  \frac{y_0}{y_1}\right) + n 
  \right\vert \right]. 
\]
Since $\Vert n\Vert_P \leq 1$ for all $n$,  
we have 
$0 \leq \log^+\left\Vert \frac{y_0}{y_1} + n \right\Vert_{P} 
\leq 
\log^+\left\Vert \frac{y_0}{y_1}\right\Vert_{P}.$
Also there is a constant $c'$ (depending on $y$) 
such that 
$\left\vert \log^+ \left\vert \sigma\left(
\frac{y_0}{y_1}\right) + n \right\vert 
- \log |n|  \right\vert \leq c'$ 
for all nonzero $n \in \ZZ$ and $\sigma$. Then there is a constant 
$c^{\prime\prime}$ (depending on $y$) such that 
$\left\vert h_{nv}(g^n(y)) - \log |n| \right\vert \leq c^{\prime\prime}$ 
for all nonzero $n \in \ZZ$. 
Hence $N_{g,y}(T) \gg\ll \exp(T)$ as $T \to +\infty$. 
\QED

For automorphisms of curves in general, we have the following. 

\begin{Lemma}
\label{lemma:elliptic}
Let $C$ be an geometrically irreducible curve over a number field $K$, 
and $f: C \to C$ an automorphism over $K$. Let $H$ be an ample divisor 
on $C$, and we fix a height function $h_H$ 
corresponding to $H$. Then we have the following. 
\begin{enumerate}
\item[(1)] 
Either \textup{(i)} $f$ is periodic, 
or \textup{(ii)} there are only finitely many $f$-periodic points on 
$C(\overline{K})$. 
\item[(2)] 
In the case \textup{(ii)}, let $x \in C(\overline{K})$ be 
a non $f$-periodic point. Then 
\[
  \# \{y \in O_f(x) \mid h_{H}(y) \leq T\} \gg\ll T 
  \quad\text{or}\quad
  \# \{y \in O_f(x) \mid h_{H}(y) \leq T\} \gg\ll \exp(T)
\]
as $T \to +\infty$. 
\end{enumerate}
\end{Lemma}

\Proof
Let $\varphi: \widetilde{C} \to C$ be the
normalization, and $\widetilde{f}: \widetilde{C} \to \widetilde{C}$ the
automorphism induced by $f$. The assertion for $(C, f)$ follows 
from that for $(\widetilde{C}, \widetilde{f})$. Hence we may 
assume that $C$ is smooth. 

{\bf Case~1:}\quad $g(C)=0$. This case is already treated 
in Lemma~\ref{lemma:p1} and Lemma~\ref{lemma:p2}. 

{\bf Case~2:}\quad $g(C)=1$. 
We fix any point $O$ of $C$ and regard $C$ as an elliptic curve 
with the origin $O$. As in the proof of Lemma~\ref{lemma:canonical:height},  
for some $n \in \ZZ_{>0}$, $f^n$ is a translation by 
a point $a \in C(\overline{K})$. If $a$ is a torsion point, 
then $f^{nm} = \id$ for some $m \in \ZZ_{>0}$, so that we have (i). 
If $a$ is not a torsion point, 
then $f^{nm}(x) = x + ma \neq x$ for any $x \in C(\overline{K})$ and 
for any  $m\in \ZZ_{>0}$.  
This shows that there are no $f^n$-periodic points. 
Hence there are no $f$-periodic points, so that we have (ii). 

Regarding (2), we will show that 
\begin{equation}
\label{eqn:elliptic}
  \# \{y \in O_f(x) \mid h_{H}(y) \leq T\} \gg\ll T
  \qquad(T\to+\infty) 
\end{equation}
for any non $f$-periodic point $x$. First we claim that we may 
replace $f$ by $f^n$. Indeed, suppose Eqn.\!~\eqref{eqn:elliptic} 
holds for $f^n$ and for any non $f^n$-periodic point. 
We will then have 
\begin{align}
\label{eqn:elliptic:2}
  & \# \{y \in O_f(x) \mid h_{H}(y) \leq T\} \\
  \notag
  & \qquad 
  = \sum_{i=0}^{n-1} \# \{y \in O_{f^n}(f^i(x)) \mid h_{H}(y) \leq T\} 
  \gg\ll T \qquad(T\to+\infty). 
\end{align}
Hence, replacing $f$ by $f^n$, it suffices to show 
\eqref{eqn:elliptic} when $f$ is a translation 
by $a$. Here $a$ is a non-torsion point, since we consider the case (ii). 
By Theorem~\ref{thm:height:machine}(vii), 
we may assume that $h_H$ is the N{\'e}ron-Tate height function $h_{NT}$ 
corresponding to an ample symmetric divisor. 
We set $2 B(x, y) :=
h_{NT}(x+y)^2 - h_{NT}(x)^2 - h_{NT}(y)^2$.  Since $B$ is a bilinear
form, we get
\[
  h_{NT}(f^n(x)) = h_{NT}(x + na) = 
  \sqrt{h_{NT}(x) + 2 n B(x, a) + n^2 h_{NT}(a)}. 
\]
Since $a$ is a non-torsion point, we have $h_{NT}(a) \neq 0$. 
Hence $\# \{y \in O_f(x) \mid h_{NT}(y) \leq T\} 
\gg\ll T$ as $T \to +\infty$. 

{\bf Case~3:}\quad $g(C)\geq 2$. 
Then $f^n = \id$ for some positive integer $n$. 
Thus we have (i). 
\QED

\subsection{$f$-periodic points}
\label{subsec:periodic}
Let $f: X \to X$ be a smooth projective surface automorphism 
over a number field $K \subset \CC$, which has positive 
topological entropy. 
Let $H$ be an ample line bundle on $X$, and 
$h_H$ a height function corresponding to $H$. 
A subset $S \subset X(\overline{K})$ is said to be 
{\em of bounded height} if there is a constant $c$ 
such that $h_H(s) \leq c$ for all $s \in S$. 

Notice that, by Theorem~\ref{thm:height:machine}(vii),  
the definition of a set of bounded height is independent 
of the choice of ample line bundles $H$. 
The Northcott finiteness property (Theorem~\ref{thm:height:machine}(v)) 
asserts that if $S$ is a set of bounded height then 
$\{s \in S \mid [K(s):K] \leq d \}$ is finite  
for every positive integer $d$. 

The set of all $f$-periodic points 
in $X(\overline{K})$ is infinite (cf. \cite[Th{\'e}or{\`e}me~in~p.~906]{Ca1}). 
We would like to know if it is a set of bounded height. 
If there is an $f$-periodic curve $C$ such that $C$ is 
pointwisely fixed by $f^n$ for some positive integer $n$, 
then $\{\textup{$f$-periodic points in $X(\overline{K})$}\}$ 
contains $C(\overline{K})$  
and so is not a set of bounded height. 
On the other hand, the following theorem asserts that, 
if there are no such curves, 
then the set of all $f$-periodic points is of bounded height. 
 
\begin{Theorem}
\label{thm:periodic:points}
Let $X$ be a smooth projective surface over a number field $K$, 
and $f: X \to X$ an automorphism of positive topological entropy. 
Let $E$ be the maximal $f$-invariant curve. 
\begin{enumerate}
\item[(1)]
The set 
$\,\left\{\textup{$f$-periodic points in 
$(X\setminus E)(\overline{K})$}\right\}\,$
is of bounded height. 
\item[(2)]
Suppose that no $f$-periodic curve is pointwisely fixed by $f^n$ 
for each $n$. Then the set 
$\,\left\{\textup{$f$-periodic points in 
$X(\overline{K})$}\right\}\,$
is of bounded height. 
\end{enumerate}
\end{Theorem}

\Proof
(1) As in the proof of Theorem~\ref{thm:canonical:height}(iii), 
let $h_Z$ be a height function corresponding to an effective 
divisor $Z$, where $Z$ satisfies $\Supp(Z) \subseteq \Supp(E)$ and 
$D - \varepsilon Z$ is ample for sufficiently small 
$\varepsilon >0$. Since 
$h_Z(x) \geq c_2$ for all 
$x \in (X\setminus E)(\overline{K})$ 
by Eqn.\!~\eqref{eqn:canonical:height:4}, we have 
\begin{align*}
\left\{\textup{$f$-periodic points in 
$(X\setminus E)(\overline{K})$}\right\}
& = 
\{x\in (X\setminus E)(\overline{K}) 
\mid
\widehat{h}_D(x) = 0\} \\
& \subseteq 
\{x\in (X\setminus E)(\overline{K})
\mid
(\widehat{h}_{D} - \varepsilon h_Z)(x) \leq - \varepsilon c_2\}. 
\end{align*}
Since $\widehat{h}_{D}- \varepsilon h_Z$ is a height function corresponding 
to an ample $\RR$-divisor, the latter set is a set of bounded height by 
Theorem~\ref{thm:canonical:height}(vi). 

(2) We use Lemma~\ref{lemma:elliptic}(1). The assumption of 
Theorem~\ref{thm:periodic:points}(2) says that the case (i) of 
Lemma~\ref{lemma:elliptic}(1) does not occur, so that  
$\left\{\textup{$f$-periodic points in 
$E(\overline{K})$}\right\}$ is a finite set. 
Hence, trivially, 
$\left\{\textup{$f$-periodic points in $E(\overline{K})$}\right\}$ 
is a set of bounded height. Now the assertion follows from (1).  
\QED

\subsection{Non $f$-periodic points}
\label{subsec:non:periodic}
We keep the notation in \S\ref{subsec:periodic}. In particular, 
$h_H$ denotes a height function corresponding to an ample divisor $H$ 
on $X$. 

Let $x$ be a non $f$-periodic point of $X(\overline{K})$. 
As in \S\ref{subsec:p1}, let $O_f(x) := 
\{f^n(x) \mid n \in \ZZ\}$ be the orbit of $x$ under $f$. 
We set 
\[
  N_{f,x}(T) := \# \{z \in O_f(x) \mid h_H(z) \leq T\}, 
\]
which depends on $T$, $O_f(x)$ and $h_H$. 
We remark that we count $N_{f,x}(T)$ with respect to $h_H$ and not 
$\widehat{h}_D$.  

Suppose further that $x$ does not lie on $E$. 
Then it follows from Theorem~\ref{thm:canonical:height}(v) 
and Proposition~\ref{prop:decomp}(2) that  
$\widehat{h}_{D_+}(z) > 0$ and $\widehat{h}_{D_-}(z) > 0$ for any 
$z \in O_f(x)$. 
Following \cite{SiK3}, 
we set 
\[
  \widehat{h}_{D}(O_f(x)) := 
  \frac{\log \left(\widehat{h}_{D_+}(z) 
  \widehat{h}_{D_-}(z)\right)}{\log\lambda}
\]
for $z \in O_f(x)$.  
We see that $\widehat{h}_{D}(O_f(x))$ is 
independent of the choice of $z \in O_f(x)$. 
 
\begin{Theorem}
\label{thm:non-periodic:points}
Let the notation and assumption be as in 
Theorem~\ref{thm:periodic:points}. 
We take a non $f$-periodic 
point $x \in X(\overline{K})$. 
Then we have the following. 
\begin{enumerate}
\item[(1)]
\begin{enumerate}
\item[(a)] 
The general case\textup{:}
If $x$ is not on $E$, then $N_{f,x}(T) \gg\ll \log T$ as 
$T \to +\infty$.  
\item[(b)]
The special case\textup{:}
If $x$ is on $E$, 
then either $N_{f,x}(T) \gg\ll T $ or $N_{f,x}(T) \gg\ll \exp(T)$ 
as $T \to +\infty$. 
\end{enumerate}
\item[(2)]
In the case \textup{(a)}, 
we have a refined estimate in terms of 
$h_{top}(f)$ and $\widehat{h}_D$\textup{:} 
\[
  N_{f,x}(T) =  \frac{2}{h_{top}(f)}\log T - \widehat{h}_D(O_f(x)) + O(1) 
  \qquad (T \to \infty),
\]
where the $O(1)$ constant is independent of $x$ 
\textup{(}but depends on $h_H$\textup{)}. 
\end{enumerate}
\end{Theorem}

\Proof
(1)(a) 
This follows from (2) below. 

(1)(b)
Let $C_0$ be an irreducible $f$-periodic curve such 
that $x \in C_0(\overline{K})$. Let $k$ be the smallest positive 
integer such that $f^k(C_0) = C_0$. We set 
$C_i := f^i(C_1)$ for $i = 1, \ldots, k-1$.  
We also set $x_0 := x$ and $x_i :=  f^i(x)$ 
for $i = 1, \ldots, k-1$.  
Since $x$ is not $f$-periodic, as in \eqref{eqn:elliptic:2}  
we get 
\[
  N_{f,x}(T) = \sum_{i=0}^{k-1} N_{f^k, x_i}(T). 
\]
By Theorem~\ref{thm:height:machine}(vii), 
there are positive constants $a_1, a_2$ and 
constant $b_1, b_2$ such that 
$a_1 h_H + b_1 \leq h_{nv} \leq a_2 h_H + b_2$ 
on every $C_i$. Since $f^k(C_i) = C_i$ and $x_i$ is not $f$-periodic, 
Lemma~\ref{lemma:elliptic} asserts 
that $N_{f^{k}, x_i}(T) \gg\ll T$ or 
$N_{f^{k}, x_i}(T) \gg\ll \exp(T)$ as $T \to +\infty$. 
Thus we get the assertion. 

(2) 
We can prove (2) as in \cite[Theorem~1.3]{SiK3}, 
together with comparison of $h_H$ and $\widehat{h}_D$. 
Indeed, we set 
\[
  \Sigma(T) := \# \{
  z \in O_f(x) \mid \widehat{h}_D(z) \leq T\}.
\] 
Since $x$ is non $f$-periodic, we have 
\begin{align*}
  \Sigma(T) 
  & = \# \{ n \in \ZZ \mid \widehat{h}_D(f^n(x)) \leq T\} \\
  & = \# \{ n \in \ZZ \mid 
  \lambda^n \widehat{h}_{D_+}(x) + 
  \lambda^{-n} \widehat{h}_{D_-}(x)\leq T\}.  
\end{align*}
It follows from \cite[Lemma~in~ p.~366]{SiK3} that, 
if $T^2 \geq 4  \widehat{h}_{D_+}(x)\widehat{h}_{D_-}(x)$, 
then 
\[
  -1 \leq 
  \Sigma(T) 
  - \frac{\log\frac{T^2}{4 \widehat{h}_{D_+}(x)\widehat{h}_{D_-}(x)}}
  {\log\lambda}
  \leq 1 + \frac{\log 4}{\log\lambda}. 
\]
Since $h_{\topo}(f) = \log\lambda$ by Theorem~\ref{thm:lambda}(1), 
we get
\begin{equation}
\label{eqn:non-periodic:points}
  \Sigma(T) = 
  \frac{2}{h_{top}(f)}\log T - \widehat{h}_D(O_f(x)) + O(1) 
  \qquad (T \to+\infty). 
\end{equation}
Since $h_H$ is a height function corresponding to 
an ample $\RR$-divisor, 
there is a positive constant $a_1 > 0$ and a constant $b_1$ 
such that
\begin{equation}
\label{eqn:non-periodic:points:2}
  h_H(x) \geq a_1 \widehat{h}_{D}(x) + b_1
  \qquad\text{for all $x \in X(\overline{K})$}.
\end{equation}
On the other hand, 
by Proposition~\ref{prop:nu:2}(2) and Eqn.\!~\eqref{eqn:canonical:height:4}, 
there is a positive constant $a_2 > 0$ and a constant $b_2$ 
such that 
\begin{equation}
\label{eqn:non-periodic:points:3}
  h_H(x) \leq a_2 \widehat{h}_{D}(x) + b_2
  \qquad\text{for all $x \in (X\setminus E)(\overline{K})$}. 
\end{equation}
Since $x$ is not on $E$, we have $O_f(x) \subset 
(X\setminus E)(\overline{K})$. 
Now the assertion follows from 
Eqns.\!~\eqref{eqn:non-periodic:points}, 
\eqref{eqn:non-periodic:points:2} and 
\eqref{eqn:non-periodic:points:3}.  
\QED

For $x \in X(\overline{K})$, let 
$O_f^+(x) := \{f^n(x) \mid n \in \ZZ_{\geq 0}\}$ denote 
the set of forward-orbit of $x$ under $f$. If 
$x$ is not $f$-periodic, we set 
$N^+_{f,x}(T) := \# \{z \in O^+_f(x) 
\mid h_H(z) \leq T\}$. Then we have the following corollary. 

\begin{Corollary}
\label{cor:non-periodic:points}
If $x \in (X\setminus E)(\overline{K})$ is not $f$-periodic, 
then 
\[
  \lim_{n \to +\infty} \frac{N^+_{f,x}(T)}{\log T} 
  = \frac{1}{h_{top}(f)}.
\] 
\end{Corollary}

\Proof
For $T \gg 1$, we set 
\[
  S(T) := \{ n \in \ZZ_{\geq 0} 
  \mid \lambda^n \widehat{h}_{D_+}(z) 
  + \lambda^{-n} \widehat{h}_{D_-}(z)\leq T\}. 
\]
If $n \in S(T)$, then $\lambda^n \widehat{h}_{D_+}(z) 
\leq T$. Thus $\# S(T) \leq 2 
+ \frac{\log T - \log \widehat{h}_{D_+}(z)}{\log\lambda}$. 
On the other hand, if $n \in \ZZ_{\geq 0}$ satisfies 
$\lambda^n \widehat{h}_{D_+}(z) + \widehat{h}_{D_-}(z)\leq T$, 
then $n \in S(T)$. Thus
$\# S(T) \geq  
\frac{\log(T - \widehat{h}_{D_-}(z)) 
- \log\widehat{h}_{D_+}(z)}{\log\lambda}$. 
By Eqns.\!~\eqref{eqn:non-periodic:points:2} 
and \eqref{eqn:non-periodic:points:3}, we get the assertion. 
\QED

\end{document}
